\documentclass[10pt]{article}
\date{}
\usepackage{hyperref}
\usepackage{amsfonts,amssymb,amsmath,amscd,euscript,array,mathrsfs}
\usepackage{enumerate}
\input{liemacs10.sty}
\newcommand{\cY}{\mathscr Y} 

\newcommand{\Part}{\mathop{{\rm Part}}\nolimits}

\renewcommand{\oline}{\overline} 
\renewcommand{\cG}{\mathcal G} 
\newcommand{\cJ}{\mathcal J} 

\newcommand{\bD}{\mathbb D}

\renewcommand{\phi}{\varphi} 
\renewcommand{\subeq}{\subseteq} 
\renewcommand{\la}{\langle}
\renewcommand{\ra}{\rangle}

\renewcommand{\mlabel}{\label}

\newcommand\ball{\mathop{\tt ball}\nolimits}

\newcommand{\Fact}{\mathop{{\rm Fact}}\nolimits}
\newcommand\Ind{\mathop{\rm Ind}\nolimits}

\begin{document} 

\title{Polynomial representations of $C^*$-algebras\\ and their applications} 
\author{Daniel Belti\c t\u a\footnote{Institute of Mathematics ``Simion Stoilow'' of the Romanian Academy, 
P.O. Box 1-764, Bucharest, Romania. Email: {\sf beltita@gmail.com, Daniel.Beltita@imar.ro}}\ 
\ and Karl-Hermann Neeb\footnote{Department of Mathematics, Friedrich-Alexander University, Erlangen-Nuremberg,
Cauerstrasse 11, 91058 Erlangen, Germany. Email: {\sf neeb@mi.uni-erlangen.de}}} 



\maketitle 


\begin{abstract} 
This is a sequel to our paper 
on nonlinear completely positive maps and dilation theory for real 
involutive algebras, where we have reduced all representation classification 
problems to the passage from a $C^*$-algebra $\cA$ to its symmetric 
powers $S^n(\cA)$, resp., to  holomorphic representations 
of the multiplicative $*$-semigroup $(\cA,\cdot)$. 
Here we study the correspondence 
between representations of $\cA$ and of $S^n(\cA)$ in detail. 
As $S^n(\cA)$ is the fixed point algebra for the natural action 
of the symmetric group $S_n$ on $\cA^{\otimes n}$, this is done 
by relating representations of $S^n(\cA)$ to those of the crossed product 
$\cA^{\otimes n} \rtimes S_n$ in which it is a hereditary subalgebra. 
For $C^*$-algebras of type I, we obtain a rather complete 
description of the equivalence classes of the irreducible 
representations of $S^n(\cA)$ and we relate this to 
the Schur--Weyl theory for $C^*$-algebras. 
Finally we show that if $\cA\subseteq B(\cH)$ is a factor of type 
II or III, then its corresponding multiplicative representation 
on $\cH^{\otimes n}$ is a factor representation of the same type, 
unlike the classical case $\cA=B(\cH)$. 
 \\ 
\textit{Mathematics Subject Classification 2010:} 
46L06, 22E66, 46L45 
\\
\textit{Keywords and phrases:} 
$C^*$-algebra,  $*$-semigroup, completely positive map
\end{abstract} 

\section{Introduction} \mlabel{sec:intro}

A real {\it seminormed involutive algebra (or $*$-algebra)} 
is a pair $(\cA,p)$, consisting 
of a real associative algebra $\cA$ endowed with an involutive antiautomorphism 
$*$ and a submultiplicative seminorm $p$ satisfying $p(a^*) =p(a)$ for $a\in \cA$. 
Our project, started in \cite{BN12}, 
eventually aims at a systematic understanding of unitary 
representations of unitary groups $\U(\cA)$ of real seminormed 
involutive algebras $(\cA,p)$. 
If $\cA$ is unital, its {\it unitary group} is 
\[ \U(\cA) := \{ a \in \cA : a^* a = aa^* = \1\}, \] 
and if $\cA$ is non-unital, then $\U(\cA):=\U(\cA^1)\cap(\1+\cA)$, 
where $\cA^1 = \cA \oplus \R \1$ is the unitization of~$\cA$.   
Typical examples we have in mind are algebras of the form 
$\cA=C^\infty(X,\cB)$ for a Banach $*$-algebra $\cB$ or  
$\cA = C^\infty(X,M_n(\K))$, where $X$ is a smooth manifold and 
$\K \in \{\R,\C,\H\}$ (or more general algebras of sections of $*$-algebra 
bundles). Focusing on unitary representations of $\U(\cA)$ 
which are boundary values of representations $\pi$ of the ball semigroups 
\[ \ball(\cA) := \{ a \in \cA \: p(a) < 1\} \]   lead us in \cite{BN15} 
to a complete reduction of this problem to the case 
where $\cA$ is a $C^*$-algebra and the representation 
is holomorphic on $\ball(\cA)$. This reduction is achieved by 
the fact that $\pi$ factors through the universal map 
$\eta_\cA \: \cA \to C^*(\cA,p)$, where 
$C^*(\cA,p)$ is the enveloping $C^*$-algebra of~$(\cA,p)$. 
\begin{footnote}{
As we learned from Jan Stochel, he obtained closely 
related results in \cite{St92} which has some overlap with \cite{BN15}.}
\end{footnote}

For a $C^*$-algebra $\cA$, we write $e^{\cA}$ 
for the $c_0$-direct sum of the $C^*$-algebras 
$S^n(\cA) \subeq \cA^{\otimes n}$, where the tensor products are constructed from 
the maximal $C^*$-cross norm. We have a holomorphic 
$*$-homomorphism of semigroups 
\[ \Gamma \: \ball(\cA) \to e^{\cA}, \quad \Gamma(a) = \sum_{n=0}^\infty a^{\otimes n}.\] 
One of the main results of \cite{BN15} 
implies that, for every bounded holomorphic $*$-re\-pre\-sen\-tation $(\pi,V)$ 
of the multiplicative $*$-semigroup $\ball(\cA)$, there exists a 
(linear) representation \break  ${\Phi \: e^\cA\to B(V)}$ 
with $\Phi \circ \Gamma = \pi$. As $e^\cA$ is the direct sum of the 
ideals $S^n(\cA)$, this in turn reduces all classification 
issues to representation theory of the $C^*$-algebras~$S^n(\cA)$. 
Via the homogeneous multiplicative map 
\begin{equation}
  \label{intro_eq1}
\Gamma_n \: \cA \to S^n(\cA), \qquad a \mapsto a^{\otimes n}, 
\end{equation}
the representations of $S^n(\cA)$ are in one-to-one correspondence 
with multiplicative holomorphic $*$-representations 
of the multiplicative $*$-semigroup $(\cA,\cdot)$ (or the unit ball $\ball(\cA)$) 
which are homogeneous of degree $n$. We call them 
{\it polynomial representations of $\cA$}. As $\U(\cA)$ is a totally
real submanifold of $\cA$ (resp., of $\1 + \cA$ in the non-unital case), 
these representations are uniquely determined by their restrictions to 
$\U(\cA)$, which are norm-continuous unitary representations. 
The multiplicative maps \eqref{intro_eq1} 
lead to natural maps 
\[ S^n(\cA)\,\hat{}\ \ssmapright{\Gamma_n^*}  \U(\cA)\,\hat{}, \qquad 
\Gamma_n^*(\pi)(a) := \pi(a^{\otimes n}), \] 
which embeds the $C^*$-algebra dual space $S^n(\cA)\,\hat{}\ $ 
into the set $\U(\cA)\,\hat{}\ $ of 
equivalence classes of unitary irreducible representations of $\U(\cA)$. 
This method of constructing points in $\U(\cA)\,\hat{}\ $ requires 
a detailed understanding of the representation theory of the 
$C^*$-algebra $S^n(\cA)$ in terms of the representation theory of $\cA$, 
and this is the main goal of the present paper. 

Representations of unitary groups are thus related to the 
interaction between $C^*$-algebras and infinite dimensional analyticity   
---a theme that was developed in the comprehensive monograph \cite{Up85} 
and has been explored until recently, for instance in the 
prequel \cite{BN15} to the present paper. 
We touch on that theme here again, 
proving, among other things, that holomorphy, although much weaker than linearity,
is strong enough to imply for instance complete positivity of 
multiplicative  homomorphisms of $C^*$-algebras 
(see Proposition~\ref{CPC} below). 
The present approach to representation theory of infinite dimensional unitary groups 
is based on $C^*$-algebras and their multiplicative  representations 
and is thus complementary to other recent works on representation theory of infinite dimensional Lie groups, 
such as \cite{Wo14} or \cite{DwOl15}. 
This approach allows us to shed fresh light on, and to partially extend,  
some results on Schur--Weyl duality in infinite dimensions from \cite{BN12}, \cite{Nes13}, and \cite{EnIz15}.  

If $\cA$ is a unital $C^*$-algebra, its unitary group $\U(\cA)$ 
is a Banach--Lie group whose representation theory is much more complex 
than the representation theory of $\cA$. 
That fact is already apparent in the special case of matrix algebras $\cA=M_k(\C)$, 
when $\U(\cA)$ is the compact Lie group $\U_k(\C)$ with a rich representation theory, while the tautological representation on $\C^k$ is the only irreducible representation of $M_k(\C)$ up to unitary equivalence. 
This classical situation is of course well understood using various well-known tools, 
including for instance 
\begin{itemize}
\item the Peter--Weyl theorem on compact groups; 
\item the Cartan--Weyl Theorem on the highest weights of irreducible representations; 
\item the Schur--Weyl theory on tensor realizations of representations of classical groups. 
\end{itemize}
These tools are quite specific to finite dimensions, inasmuch as they rely on Haar measures 
or on eigenvectors of the adjoint representation. 

We have shown in \cite{BN12} that,  
for any unital $C^*$-algebra $\cA$, the Schur--Weyl method leads 
for every irreducible representation $(\pi,\cH)$ of $\cA$ 
to an infinite family $(\pi_\lambda, \bS_\lambda(\cH))$ 
of irreducible representations of $\U(\cA)$, where 
$\lambda \in \Part(n)$ is a partition of $n$ and 
$\cH^{\otimes n}$ decomposes under $\U(\cA)$ 
as $\bigoplus_{\lambda \in \Part(n)} \bS_\lambda(\cH)^{\oplus d_\lambda}$, 
where $d_\lambda$ is the dimension of the irreducible representation 
of $S_n$ corresponding to~$\lambda$. 
In \cite{BN12} it was left unclear to which 
extent one thus exhausts the unitary dual $\U(\cA)\,\hat{}\ $ 
of $\U(\cA)$. 
As a by-product of our investigation, we prove in Theorem~\ref{SW-main} that if $\cA$ is a separable $C^*$-algebra, 
then the unitary dual $\U(\cA)\,\hat{}\ $ 
is exhausted by the representations constructed by the Schur--Weyl method if and only if $\cA$ is the $C^*$-algebra $K(\cH)$ 
 of compact operators on a separable complex Hilbert space. 
Among other things, this explains why the remarkably complete results of \cite{Ki73} on unitary irreducible representations 
of $\U(\cA)$ for $\cA=K(\cH)$ 
were never extended to unitary groups of more general $C^*$-algebras. 

Our approach to the 
representation theory of $S^n(\cA)$ is based on the fact that this algebra is the fixed-point set of the permutation action of the permutation group $S_n$ on the tensor power $\cA^{\otimes n}$, 
hence we can use the old observation of \cite{Ro79} for describing the dual space $S^n(\cA)\,\hat{}\ $ 
as an open subset of the dual space of the corresponding crossed product $\cA^{\otimes n}\rtimes S_n$. 
We actually extend that approach to unitary equivalence classes of 
factor representations, 
because the unitary group $\U(\cA)$ in general has norm-continuous 
unitary factor representations of type II and III. 
In this connection, some specific applications are given in Theorem~\ref{fact3}. 

The structure of this paper is as follows. 
In Section~\ref{sec:2} we discuss the structure of the symmetric powers 
$S^n(\cA)$ of a few specific types of $C^*$-algebras. 
In Section~\ref{sec:3} we then take a closer look at 
multiplicative homomorphisms, resp., representations 
of $C^*$-algebras. Here our main result is Theorem~\ref{CPC},  
asserting that every multiplicative homomorphism of 
$C^*$-algebras is completely positive, contractive and a
sum of homogeneous homomorphisms. 
To understand the passage from $\cA$ to $S^n(\cA)$, it is natural 
to study  $S^n(\cA)$ as the fixed point algebra $(\cA^{\otimes n})^{S_n}$ 
for the natural action of the symmetric group $S_n$ on~$\cA^{\otimes n}$. 
Therefore we 
study in Section~\ref{sec:4} $C^*$-dynamical systems 
$(\cA,G,\alpha)$, where $G$ is a finite group. Here our main point 
is to describe the relation between representations of the fixed point 
algebra $\cA^G$ and the crossed product $\cA \rtimes_\alpha G$, in which 
$\cA^G$ is embedded as a hereditary subalgebra. 
Here the main results are the description of the 
factor representations of the crossed product 
in terms of induced representations (Theorem~\ref{factbij}), 
its refinement for irreducible representations 
(Proposition~\ref{irr}) and, eventually, the description 
of $\hat{\cA^G}$ in Theorem~\ref{fix}.
All this is applied in Section~\ref{sec:5} to the 
special case $(\cA^{\otimes n}, S_n, \alpha)$. 
For $C^*$-algebras of type I, we obtain in Theorem~\ref{exact} 
a description of the equivalence classes of the irreducible 
representations of $S^n(\cA)$ in terms of $\cA^{\otimes n} \rtimes_\alpha S_n$ 
and we also relate this to the Schur--Weyl theory for $C^*$-algebras 
developed in \cite{BN12}. 
In the final Section~\ref{sec:6} we develop some 
aspects of Schur--Weyl theory for factor representations. 
This generalizes parts of recent work of Nessonov \cite{Nes13} 
and Enomoto/Izumi \cite[Th.~4.1]{EnIz15} 
to more general von Neumann algebras. 
Our main results (Proposition~\ref{fact2} and Theorem~\ref{fact3}) 
imply that the multiplicative $n$-fold tensor power representation on $\cH^{\otimes n}$  for a factor $\cM \subeq B(\cH)$ of type II$_1$, II$_\infty$  or III 
is of the same type and we even obtain some more detailed information 
on the corresponding projections.

{\bf Basic notation.} 
We will use the following notation for any $C^*$-algebra~$\cA$: 
\begin{itemize}
\item $\U(\cA):=\{u\in\cA\mid u^*u=uu^*=\1\}$ for the \emph{unitary group} of a unital 
$C^*$-al\-ge\-bra $\cA$ and for a non-unital $C^*$-algebra we put $\U(\cA):=\U(\cA^1)\cap(\1+\cA)$, 
where $\cA^1 = \cA \oplus \C \1\subseteq M(\cA)$ is the unitization of $\cA$, 
and $M(\cA)$ is the multiplier algebra of~$\cA$; 
\item $\ball(\cA):=\{a\in\cA\mid \Vert a\Vert< 1\}$ for the \emph{open ball semigroup} of $\cA$; 
\item $\oline{\ball}(\cA):=\{a\in\cA\mid \Vert a\Vert \le1\}$ for the \emph{closed ball semigroup} of $(\cA,p)$.   
\end{itemize}
Both $\ball(\cA)$ and $\oline{\ball}(\cA)$ are $*$-semigroups, 
and if $\cA$ is unital, then $\1 \in\oline{\ball}(\cA,p)$. 
If $\cH$ is a complex Hilbert space and $B(\cH)$ is the $C^*$-algebra of all 
bounded linear operators on $\cH$, 
then we also write $C(\cH):=\oline{\ball}(B(\cH))=\{T\in B(\cH)\mid \Vert T\Vert\le 1\}$.

We denote by $\simeq$ the relation of unitary equivalence and by $\approx$ the (weaker) relation of quasi-equivalence 
of group or algebra representations on Hilbert spaces \cite{Dix64}. 
For any $C^*$-algebra $\cA$ we denote by $\widehat{\cA}$ its set of unitary equivalence classes of irreducible $*$-representations. 
For every irreducible $*$-representation $\pi\colon\cA\to B(\cH_\pi)$, 
we denote its unitary equivalence class by 
$[\pi]\in\widehat{\cA}$. 
The set $\widehat{\cA}$ is considered as a topological space endowed with the Fell topology, 
which may be described by the condition that   
the correspondence $\cJ\mapsto\{[\pi]\in\widehat{\cA}\mid \cJ\not\subset\ker\pi\}$ is a bijection between the closed two-sided ideals $\cJ\subseteq\cA$ and the open subsets of $\widehat{\cA}$. 

\tableofcontents

\section{Some examples of symmetric powers 
of $C^*$-algebras}  \mlabel{sec:2}

As mentioned above, in the present paper we develop some mechanisms to get more explicit 
access to the $C^*$-algebras $S^n(\cA)$ for a given $C^*$-algebra~$\cA$. 
To illustrate that idea, we will briefly discuss here symmetric tensor powers 
of some basic examples of $C^*$-algebras.

\begin{ex}\mlabel{ex1}
If $X$ is a compact space and $\cA = C(X)$, then, 
for every $n \in \N$, we have  
$\cA^{\otimes n} \cong C(X^n)$ because 
$\Hom(\cA^{\otimes n},\C) \cong \Hom(\cA,\C)^n$. 
From that we immediately infer that 
\[ S^n(\cA) \cong C(X^n/S_n)\] 
where $X^n/S_n$ stands for the quotient space with respect to the natural action of the symmetric group $S_n$ 
by permutations of the coordinates of points in the $n$th Cartesian power~$X^n$. 
Note that $X^n/S_n$ need not be a smooth manifold if $X$ is. 

If $X=\{x_1,\dots,x_N\}$ is a finite set, then $C(X)\cong \C^N$, 
hence the above remarks imply that for any $n\ge1$, we have a linear isomorphism 
\[S^n(\C^N)\cong \C^{r(N,n)}.\] 
Here $r(N,n):=\vert X^n/S_n\vert=\binom{N+n-1}{n}=\frac{(N+n-1)!}{n!(N-1)!}$ is the number of combinations of $N$ elements taken by $n$ at a time, 
where we allow any element to be selected repeatedly 
(equivalently, the number of non-decreasing functions $\{1,\dots,n\}\to\{1,\dots,N\}$). 
\end{ex}

\begin{ex}\mlabel{ex2} 
Let $\cA := C(X,\cB)$ for a compact space $X$ and a unital 
$C^*$-algebra $\cB$. 
It is well known that $C(X)$ is nuclear and one has a canonical $*$-isomorphism 
$\cA \cong C(X) \otimes \cB$ 
(\cite[Th. II.9.4.4]{Bl06}, \cite[Thm.~6.4.17]{Mu90}). 
We 
want to determine the $C^*$-algebra $S^n(\cA)$. It is easy to see that 
\[ \cA^{\otimes n} \cong (C(X) \otimes \cB)^{\otimes n} 
\cong C(X^n, \cB^{\otimes n}),\] 
where $S_n$ acts on this algebra by 
\[ (\sigma F)(x_1, \ldots, x_n) := \sigma.F(x_{\sigma^{-1}(1)}, \ldots, x_{\sigma^{-1}(n)}).\]
Accordingly, $S^n(\cA)$ can be identified with the set 
$C(X^n, \cB^{\otimes n})^{S_n}$ of 
$S_n$-equivariant maps $F \: X^n \to \cB^{\otimes n}$. 
In this sense it is an {\it equivariant map algebra} (cf.\ \cite{NS12}). 

The center of $\cA^{\otimes n}$ contains $C(X^n)$, so that we obtain 
\[ Z(S^n(\cA)) \supeq  S^n(C(X)) \cong C(X^n)^{S_n} \cong C(X^n/S_n)\] 
(see \cite[Prop.~12.2]{St92} for more details). 
Since every irreducible representation of the $C^*$-algebra $S^n(\cA)$ determines a central 
character $\chi \: S^n(C(X)) \to \C$, and all these characters are given by 
evaluation in an $S_n$-orbit 
$\bx := S_n.(x_1,\ldots, x_n) =: [x_1, \ldots, x_n] \in X^n/S_n$, it factors 
through the restriction map given by 
\[ \ev_{\bx} \: S^n(\cA) \cong C(X^n, \cB^{\otimes n})^{S_n} 
\to C(\bx, \cB^{\otimes n})^{S_n}, \quad 
f \mapsto f\res_{\bx}.\] 
Let $S_{n,\bx}$ be the stabilizer of some representative 
$(x_1,\ldots, x_n) \in \bx$, so that $\bx \cong S_n/S_{n,\bx}$. 
Evaluating in the representative, 
we see that the image of $\ev_{\bx}$ is isomorphic to 
$(\cB^{\otimes n})^{S_{n,\bx}}$. 
Assuming without loss of generality that $|\{x_1, \ldots, x_n\}| = k$, 
we see that 
\[ S_{n,\bx} \cong S_{n_1} \times \cdots \times S_{n_k} \quad \mbox{ 
for } \quad n_1 + \cdots + n_k = n \] 
is a product of symmetric groups. 
This further leads to 
\[ C(\bx, \cB^{\otimes n})^{S_n} \cong (\cB^{\otimes n})^{S_{n,\bx}} 
\cong \bigotimes_{j = 1}^k (\cB^{\otimes n_j})^{S_{n_j}}
\cong \bigotimes_{j = 1}^k S^{n_j}(\cB).\] 
Therefore the determination of the irreducible representations 
of $S^n(\cA)$ is reduced to the determination of irreducible representations 
of finite tensor products of the algebras $S^m(\cB)$. 
If $\cB$ is of type $I$, then 
so is every $S^m(\cB)$ (see Proposition~\ref{cross}\eqref{cross_item5} below). 
Further, irreducible representations of tensor powers of type I algebras 
are tensor products of irreducible representations. 
This reduces the problem 
to determine the irreducible representations of $S^n(\cA)$ to the corresponding 
problem for the target algebras $S^n(\cB)$. 
\end{ex}

From \cite{BN12} we know that every irreducible representation 
$(\pi, \cH)$ of a $C^*$-algebra $\cA$ 
and every partition $\lambda$ of $n$ determine an irreducible representation 
$(\pi_\lambda, \bS_\lambda(\cH))$ of $S^n(\cA)$. 
One may expect to obtain from every 
faithful irreducible representation $(\pi, \cH)$ of $\cA$ a decomposition 
\begin{equation}\label{hope}
S^n(\cA) \cong \bigoplus_{\lambda \in \Part(n)} S^n_\lambda(\cA)
\end{equation}
from the corresponding decomposition of $S^n(B(\cH))$, acting on $\cH^{\otimes n}$, 
where it commutes with the action of $S_n$. 
We shall see in Example~\ref{ex:bh} below that this is 
far from being true in general, as the example 
$\cA = B(\cH)$ shows. However, it is true for the ideal 
$K(\cH)$ of compact operators.

\begin{ex}\mlabel{ex3}
Let $\cA=K(\cH)$ be the $C^*$-algebra of compact operators on the infinite 
dimensional complex Hilbert space~$\cH$. 
Using the fact that for any finite dimensional complex vector spaces $\cV_1$ and $\cV_2$ one has a canonical isomorphism 
$\End(\cV_1)\otimes\End(\cV_2)\cong \End(\cV_1\otimes\cV_2)$ 
and approximating the compact operators by finite-rank operators, 
we obtain, 
for any integer $n\ge 1$, a canonical isomorphism of 
$C^*$-algebras $\cA^{\otimes n}\cong K(\cH^{\otimes n})$ (see \cite[Ex.~B.20]{RW98} 
for details). This isomorphism 
intertwines the permutation action of $S_n$ on $\cA^{\otimes n}$ 
and the spatial action of $S_n$ on $K(\cH^{\otimes n})$ defined via the unitary representation 
\[ \tau\colon S_n\to \U(\cH^{\otimes n}), \quad 
\tau(\sigma)(v_1 \otimes \cdots \otimes v_n) 
= v_{\sigma^{-1}(1)} \otimes \cdots \otimes v_{\sigma^{-1}(n)}. \] 
defined by permutation of factors in tensor monomials. 
Therefore, considering the fixed-point algebras of the above two actions of $S_n$, 
we see that the isomorphism $\cA^{\otimes n}\cong K(\cH^{\otimes n})$ defines by restriction an isomorphism 
\begin{equation}
  \label{eq:sn-deco}
S^n(\cA)\cong K(\cH^{\otimes n})\cap\tau(S_n)' \cong \bigoplus_{\lambda \in \Part(n)} 
K(\bS_\lambda(\cH)),
\end{equation}
where $\bS_\lambda(\cH) = B(V_\lambda, \cH^{\otimes n})^{S_n}$ is the multiplicity 
space for the irreducible representation $(\rho_\lambda, V_\lambda)$ of $S_n$ 
in $\cH^{\otimes n}$. In fact, from the $S_n$-decomposition 
\[ \cH^{\otimes n} \cong \bigoplus_{\lambda \in \Part(n)} 
\bS_\lambda(\cH)\otimes V_\lambda \] 
it follows that a compact operator on $\cH^{\otimes n}$ commutes with 
$S_n$ if and only if it is a direct sum of operators of the form 
$A_\lambda \otimes \1$, where $A_\lambda \in B(\bS_\lambda(\cH))$. 
Note that \eqref{eq:sn-deco} is the 
decomposition of $S^n(\cA)$ according to the general 
structure theory of $C^*$-algebras of compact operators 
as developed for instance in \cite[Thm.~1.4.5]{Ar76}. 
\end{ex}

\begin{ex}\mlabel{ex4}
If $\cA$ and $\cB$ are any 
 $C^*$-algebras and $\phi\colon\cA\to\cB$ is a $*$-homomorphism, 
then for any integer $n\ge 1$ one has natural $*$-homomorphisms 
$$\phi^{\otimes n}\colon\cA^{\otimes n}\to\cB^{\otimes n}\quad\text{ and } \quad 
S^n(\phi)\colon S^n(\cA)\to S^n(\cB),$$ 
and $\bullet^{\otimes n}$ and $S^n(\bullet)$ are functors 
from the category of $C^*$-algebras into itself. 
The behavior of these functors with respect to inductive limits of nuclear $C^*$-algebras  
sheds some light on symmetric tensor powers of AF-algebras, 
that is, $C^*$-algebras that are inductive limits of finite dimensional $C^*$-algebras. 

In fact, let $\cA$ be any nuclear $C^*$-algebra with a sequence of 
nuclear closed $*$-subalgebras $\cA_1\subseteq\cA_2\subseteq\cdots\subseteq\cA$ 
whose union is dense in $\cA$. 
It is easily seen that one has 
a sequence of closed $*$-subalgebras $\cA_1^{\otimes n}\subseteq\cA_2^{\otimes n}\subseteq\cdots\subseteq\cA^{\otimes n}$ 
whose union is dense in $\cA^{\otimes n}$, and then 
a sequence of closed $*$-subalgebras $S^n(\cA_1)\subseteq S^n(\cA_2)\subseteq\cdots\subseteq S^n(\cA)$ 
whose union is dense in $S^n(\cA)$. 
Thus each of the functors $\bullet^{\otimes n}$ and $S^n(\bullet)$ commutes with inductive limits of nuclear $C^*$-algebras. 
Since these functors preserve the category of finite dimensional $C^*$-algebras, it then follows that they also preserve the category of AF-algebras. 

We point out however that, as a by-product of Example~\ref{ex3}, if $\cA$ is a UHF-algebra 
(that is, a unital $C^*$-algebra that is the inductive limit of a sequence of unital $*$-subalgebras isomorphic to full matrix algebras) 
and $n\ge 1$, then $S^n(\cA)$ is not a UHF-algebra. 
More precisely, Example~\ref{ex3} shows that if $\cA=M_k(\C)$ for some $k\ge 2$, then $S^n(\cA)$ is the direct sum of at least two full matrix algebras if $n\ge 2$. 

Now let $\cA \cong \otimes_{k = 1}^\infty B(E_k)$ with $E_k \cong \C^{n_k}$ 
be an UHF algebra. Since $\cA$ is a direct limit of finite dimensional simple matrix algebras 
$\cA_N \cong B(F_N)$ with 
$F_N \cong E_1 \otimes \cdots \otimes E_N$, 
the $C^*$-algebra $S^n(\cA)$ is the direct limit of the 
$C^*$-algebras $S^n(\cA_N) \cong \bigoplus_{\lambda \in \Part(n)} 
S^n(\cA_N)_\lambda$ with 
$S^n(\cA_N)_\lambda \cong B(\bS_\lambda(F_N))$ 
(Example~\ref{ex3}), where the projection onto the $\lambda$-summand 
corresponds to restriction to the corresponding $S_n$-isotypic subspace 
in $F_N^{\otimes n}$. From the $S_n$-equivariance of the inclusions 
$F_N^{\otimes n}\into F_{N+1}^{\otimes n}$, it follows that the inclusion 
$\cA_N \into \cA_{N+1}$ induces inclusions 
$S^n(\cA_N)_\lambda \into S^n(\cA_{N+1})_\lambda$, so that we obtain a 
direct sum decomposition 
\[ S^n(\cA) \cong \bigoplus_{\lambda \in \Part(n)} S^n(\cA)_\lambda \quad \mbox{ with }\quad 
S^n(\cA)_\lambda \cong \indlim S^n(\cA_N)_\lambda, \] 
and the construction implies in particular that the algebras 
$S^n(\cA)_\lambda$ are simple and act irreducibly 
on $\bS_\lambda(F)$, where $F$ is the Hilbert direct limit of the spaces~$F_N$.
Therefore we find a decomposition similar to the one 
in Example~\ref{ex3}. 
\end{ex}

\begin{ex} \mlabel{ex:bh} 
Let $\cH$ be an infinite dimensional Hilbert space. 

(a) {\it The representation of $S^2(B(\cH))$ on $\cH^{\otimes 2}$ is 
not injective.} First we recall the non-uniqueness of $C^*$-norms 
on $B(\cH)\otimes B(\cH)$ from \cite{JP95} which implies 
that the spatial tensor product 
(corresponding to the least $C^*$-norm by \cite[Thm.~6.4.18]{Mu90}), 
which is the image $B(\cH)^{\bar{\otimes} 2}$ of $B(\cH)^{\otimes 2}$ under 
the natural representation on $\cH^{\otimes 2}$ 
differs from  the tensor product $B(\cH)^{\otimes 2}$ 
with respect to the maximal $C^*$-norm. This implies that the 
kernel $\cJ \subeq B(\cH)^{\otimes 2}$ of the representation on 
$\cH^{\otimes 2}$ is non-zero. Clearly, $\cJ$ is invariant under the 
flip automorphism $\tau(a_1 \otimes a_2) = a_2 \otimes a_1$. 
If $0 \leq a \in \cJ$ is non-zero, then the same holds for 
$a + \tau(a) \in \cJ \cap S^2(B(\cH))$. Therefore the 
representation of $S^2(B(\cH))$ on $\cH^{\otimes 2}$ is also not injective. 

(b) {\it The representation of $S^2(B(\cH))$ on $\cH^{\otimes 2}$ is 
not surjective onto $B(\cH^{\otimes 2})^{S_2}$}. 
Let $(e_j)_{j \in J}$ be an orthonormal basis of $\cH$ and 
$P \in B(\cH^{\otimes 2})$ be the orthogonal projection onto the subspace 
spanned by $(e_j \otimes e_j)_{j \in J}$. Then $P$ commutes with the flip 
action of $S_2$ 
and it remains to show that it is not contained in the spatial tensor product 
$B(\cH)^{\oline{\otimes} 2}$. To this end, we write $\cH^{\otimes 2} \cong \oplus_{j \in J} 
(\cH_j \otimes e_j)$ as a direct sum and, accordingly, $P$ as a matrix 
$(P_{ij})_{i,j \in J}$. Then $P_{ij}(v)= \delta_{ij} \la v, e_j \ra e_j$ 
is either $0$ or the orthogonal projection onto 
$\C (e_j \otimes e_j)$. In particular, the subset 
$\{ P_{ij} \: i,j \in J\} \subeq B(\cH)$ is closed, discrete and 
infinite, hence not relatively compact. 
According to Exercise~11.5.7 in  \cite{KR92}, 
this implies that $P$ is not contained in $B(\cH)^{\bar\otimes 2}$. 
\end{ex}

\section{Multiplicative homomorphisms of $C^*$-algebras}
\mlabel{sec:3}

In this section we take a closer look at 
multiplicative homomorphisms, resp., representations 
of $C^*$-algebras. Here our main result is Theorem~\ref{CPC},  
asserting that every multiplicative homomorphism of 
$C^*$-algebras is completely positive, contractive and a
sum of homogeneous homomorphisms. 
We shall use 
some results and notions from our earlier paper \cite{BN15} 
on completely positive nonlinear maps and their dilation theory. 

\begin{defn}\label{nonlin}
If $\cA$ and $\cB$ are any unital $C^*$-algebras, 
then a map $\phi\colon\cA\to\cB$ is a  \emph{multiplicative homomorphism} 
if it is holomorphic and satisfies 
\[ \phi(\1)=\1, \qquad \phi(ab)=\phi(a)\phi(b) \quad \mbox{ and } \quad 
\phi(a^*)=\phi(a)^* \quad \mbox{ for } \quad a,b\in\cA.\] 
A linear map $\psi\colon \cA\to\cB$ is called \emph{contractive} if it satisfies  $\psi(\oline{\ball}(\cA))\subseteq\oline{\ball}(\cB)$.
\end{defn}

\begin{thm}\mlabel{CPC}
Every  multiplicative homomorphism of unital $C^*$-algebras $\phi\colon\cA\to\cB$ 
is completely positive, contractive and has a weakly convergent 
series expansion $\phi=\sum\limits_{n\ge0}\phi_n$, 
where $\phi_n$ is multiplicative and homogeneous of degree~$n$.
\end{thm}

The proof of the above result is based on the following observations. 
For the case of a positively homogeneous map $\phi$, the complete positivity 
of $\phi$ can also be derived from \cite[Prop.~11.1]{St92}.

\begin{lem} \mlabel{prop:7.1spec} 
Every dilatable bounded completely positive function 
$\phi \: \ball(\cA) \to B(V)$ extends to a weakly continuous 
completely positive function $\oline\phi\colon \oline{\ball}(\cA)\to B(V)$. 
\end{lem}

\begin{prf} 
First recall from \cite{BN15} that 
$S^\circ := \ball(\cA)$ is a semigroup ideal in 
$S := \oline{\ball}(\cA)$. 
Let $(\pi, \cH,\iota)=(\pi_\phi,\cH_\phi,\iota_\phi)$ be the minimal dilation of~$\phi$ 
given by 
\cite[Prop. 3.4]{BN15}. 
Then we use  
\cite[Prop. 3.21]{BN15}
to extend the representation 
$\pi$ to a weakly continuous representation $\widehat{\pi}$ of $S$, 
where the weak continuity follows by the explicit formula of~$\widehat{\pi}$, 
since $\phi$ is norm-continuous by 
\cite[Th. 6.4]{BN15}. 
Note that the representation $\pi$ is completely positive by  
\cite[Prop. 4.1]{BN15}, 
hence, by continuity, also $\widehat{\pi}$ is completely positive. 
Then 
$\oline\phi(s) := \iota^* \widehat{\pi}(s) \iota$ is a completely positive 
weakly continuous extension of $\phi$ to $S$. 
\end{prf} 

\begin{prop}\mlabel{sCPC}
Let $\cA$ be any unital $C^*$-algebra, and 
$\pi\colon\ball(\cA)\to B(\cH_\pi)$ 
be any holomorphic $*$-representation.  
Then $\pi$ is 
completely positive, and extends to a weakly continuous $*$-representation 
$\oline\pi\colon \oline{\ball}(\cA)\to C(\cH_\pi)$. 
If $\pi$ is nondegenerate, then $\lim\limits_{r\to1-}\pi(r\1)=\1$ in the weak operator topology. 
\end{prop}

\begin{prf} 
The map $\pi$ is positive definite because it is a homomorphism of $*$-semigroups. 
Therefore, after we will have proved that $\pi$ is contractive, hence in particular bounded, 
it will follow by \cite[Th. 1.1]{BN15} 
that it is completely positive. 
Then Lemma~\ref{prop:7.1spec}  
ensures that $\pi$ extends to a weakly continuous $*$-representation 
of $\oline{\ball}(\cA)$.

The proof of the fact that $\pi$ is contractive has two stages. 

{\bf Step 1:} 
We first consider the case $\cA=\C$, hence $\ball(\cA)=\bD$ is the unit disk. 
Then there exists the norm-convergent power series expansion 
$\pi(z)=\sum_{n\ge0}z^nP_n\in B(\cH_\pi)$ for 
all $z\in\bD$, with $P_n\in B(\cH_\pi)$ for all $n\ge 0$. 
For all $z,w\in\bD$ one has $\pi(\bar z)=\pi(z)^*$ and $\pi(zw)=\pi(z)\pi(w)$, 
and this implies that the coefficients of the above power series satisfy 
the conditions $P_n=P_n^*=P_n^2$ and $P_nP_m=0$ for all $n,m\ge 0$ with $n\ne m$. 
That is, $(P_n)_{n\ge 0}$ is a sequence of mutually orthogonal projections, 
and then 
$$(\forall z\in\bD)\quad \Vert\pi(z)\Vert\le\sup_{n\ge 0}\vert z^n\vert=1$$
which in particular shows that $\pi$ is bounded. 
Moreover, one has $\lim\limits_{0<r\to1-}\pi(r\1)=\sum\limits_{n\ge 0}P_n$ in the weak operator topology in $B(\cH_\pi)$. 

{\bf Step 2:} 
In the general case, for all $z\in\bD$ and $a\in\ball(\cA)$ one has 
$\pi(za)=\pi(z\1)\pi(a)$, hence if $\pi\not\equiv0$, then necessarily $\pi\vert_{\bD}\not\equiv0$, 
which we may assume until the end of this proof. 
Then use the above stage to find the sequence of mutually 
orthogonal projections $(P_n)_{n\ge 0}$ in $B(\cH_\pi)$ 
with $\pi(z\1)=\sum\limits_{n\ge0}z^nP_n\in B(\cH_\pi)$ for 
all $z\in\bD$. 

For each $a\in\ball(\cA)$ and all $z\in\bD$ one has 
$\pi(z\1)\pi(a)=\pi(a)\pi(z\1)$, 
and this implies $P_n\pi(a)=\pi(a)P_n$, 
hence $\pi_n(\cdot):=P_n\pi(\cdot)$ is a holomorphic $*$-representation of $\ball(\cA)$.  
Moreover $\pi_n(za)=z^n\pi_n(a)$ for all $a\in\ball(\cA)$ and $z\in\bD$, 
hence one may use \cite[Lemma 2.11]{BN15} to extend $\pi_n$ 
to a holomorphic $*$-representation $\pi_n\colon(\cA,\cdot)\to B(\cH_\pi)$ 
which is homogeneous of degree~$n$. 
Since $\pi_n$ is continuous at $0\in\cA$, there exist $r,M>0$ with $\Vert\pi_n(a)\Vert\le M$ 
if $a\in\cA$ with $\Vert a\Vert\le r$. 
Since $\pi_n$ is homogeneous, it then follows that it is bounded on $\ball(\cA)$, 
and then \cite[Rem. 3.5]{BN15} implies $\Vert\pi_n(a)\Vert\le 1$ for all $a\in\ball(\cA)$.  

Denoting $P:=\sum\limits_{n\ge0}P_n$, 
we then obtain the weakly convergent series 
\[ \sum\limits_{n\ge0}\pi_n(a)=P\pi(a)=\pi(a) \quad \mbox{ and } \quad 
\Vert\pi(a)\Vert\le\sup\limits_{n\ge 0}\Vert\pi_n(a)\Vert\le 1 \] 
for every $a\in\ball(\cA)$. 
The condition that $\pi$ be nondegenerate means $P=\1$, and the assertion follows. 
\end{prf}

Finally, we can now prove Theorem~\ref{CPC}.

\begin{prf} 
Using a realization of $\cB$ as an operator algebra, we may actually assume $\cB=B(\cH)$ for some complex Hilbert space~$\cH$. 
Then Proposition~\ref{sCPC} and its proof show 
that the restricted map 
$\pi:=\phi\vert_{\ball(\cA)}\colon\ball(\cA)\to B(\cH)$ 
is completely positive and there exists a sequence of mutually orthogonal projections $(P_n)_{n\ge 0}$ 
such that $\pi_n:=\pi(\cdot)P_n=P_n\pi(\cdot)$ is homogeneous of degree~$n$ and completely positive on $\ball(\cA)$, 
hence in fact $\pi_n\colon\cA\to B(\cH)$ is completely positive. 

On the other hand, since $(P_n)_{n\ge 0}$  are mutually orthogonal projections, we have 
$\phi=\sum\limits_{n\ge0}\pi_n$, 
this series being convergent in the pointwise weak operator topology. 
Therefore also $\phi$ is completely positive. 
\end{prf}

\begin{rem}\mlabel{square}
In connection with Theorem~\ref{CPC} we point out that for some simple holomorphic maps on a $C^*$-algebra~$\cA$ 
the complete positivity property is a condition that  
has strong structural implications on~$\cA$.  
For instance, if $\cA$ is unital and we define $\phi\colon\cA\to\cA$, $\phi(a):=a^2$, 
then $\phi$ is a holomorphic positive contractive map,  
and the following conditions are equivalent: 
\begin{enumerate}[(i)]
\item\mlabel{square:i} The map $\phi$ is completely positive. 
\item\mlabel{square:ii} The map $\phi$ is a multiplicative homomorphism. 
\item\mlabel{square:iii} The $C^*$-algebra $\cA$ is commutative. 
\end{enumerate}
In fact, one has \eqref{square:i}$\Leftrightarrow$\eqref{square:iii} for instance by \cite[Th. 3]{JiT03}, 
and it is clear that \eqref{square:iii}$\Rightarrow$\eqref{square:ii}. 
Conversely, if we assume that $\phi$ is a multiplicative homomorphism, then for arbitrary $a,b\in\cA$ 
and $t,s\in\R$ one has $(e^{ta}e^{sb})^2=\phi(e^{ta}e^{sb})=\phi(e^{ta})\phi(e^{sb})=(e^{ta})^2(e^{sb})^2$, 
which is equivalent to $e^{sb}e^{ta}=e^{ta}e^{sb}$. 
Differentiating this equality with respect to $s$ and $t$ at $s=t=0$, we obtain $ba=ab$, hence the $C^*$-algebra $\cA$ is commutative. 
\end{rem}

\section{Representations of crossed products by finite groups}
\mlabel{sec:4}

In this section, unless otherwise mentioned, 
$(\cA,G,\alpha)$ is a $C^*$-dynamical system, where $G$ is a finite group 
that acts by $*$-automorphisms of a $C^*$-algebra $\cA$ by $(g,x)\mapsto \alpha_g(x)$. 
We investigate here how the covariant representations of these data 
are related to the representations of the fixed-point subalgebra 
\[ \cA^G:=\{x\in\cA\mid (\forall g\in G)\ \alpha_g(x)=x\},\] 
which by \cite{Ro79} is $*$-isomorphic 
to a hereditary subalgebra of the crossed product ${\cA\rtimes_\alpha G}$. 
These results will be later used for maximal $C^*$-tensor powers $\cA^{\otimes m}$ 
acted on by the permutation group $G=S_m$, with the aim of describing irreducible 
or factor representations of the symmetric $C^*$-tensor powers $(\cA^{\otimes m})^G
=S^m(\cA)$, for any $C^*$-algebra~$\cA$.

The \emph{crossed product} corresponding to the above $C^*$-dynamical system is the $C^*$-algebra 
$\cA\rtimes_\alpha G:=\ell^1(G,\cA;\alpha):=\{f\mid f\colon G\to\cA\}$ 
with the multiplication 
$$(f_1f_2)(g):=\frac{1}{\vert G\vert}\sum\limits_{h\in G} f_1(h)\alpha_h(f_2(h^{-1}g))$$ 
and the involution $f^*(g):=\alpha_g(f(g^{-1}))^*$ for all $f,f_1,f_2\in \ell^1(G,\cA;\alpha)$ 
and $g\in G$. 
A \emph{covariant representation} is a pair $(\pi,U)$, where  
$\pi\colon \cA \to B(\cH)$ is a $*$-representation and  
$U\colon G \to B(\cH)$ is a unitary representation satisfying 
$\pi(\alpha_g(x))=U_g \pi(x)U_g^*$ for all $x\in\cA$ and $g\in G$.  
We also define the $*$-representation 
$$\pi{\rtimes U}\colon \cA\rtimes_\alpha G\to B(\cH),\quad 
(\pi{\rtimes U})(f)=\sum\limits_{g\in G}\pi(f(g))U_g.
$$
We say that the covariant representation $(\pi,U)$ is a \emph{factor representation} 
if the von Neumann algebra generated by $(\pi{\rtimes U})(\cA\rtimes_\alpha G)$, 
that is $(\pi(\cA)\cup\pi(G))''$, is a factor etc.

The following definition is a specialization of \cite[Def. 3.2]{Tak67} to finite groups, 
except that we use here the left regular representation of groups. 

\begin{defn}\mlabel{Def3.1} (induced representations)
Let $(\cA,G,\alpha)$ be as above. 
Assume that $G_0$ is a subgroup of the finite group $G$, 
denote by $(\cA,G_0,\alpha_0)$ the corresponding $C^*$-dynamical system with $\alpha_0:=\alpha\vert_{G_0}$. 

Then for any covariant representation $(\pi_0,V)$ of $(\cA,G_0,\alpha_0)$ on some Hilbert space~$\cH_0$, 
its corresponding \emph{induced representation} is the covariant representation $(\pi,U)$ 
of $(\cA,G,\alpha)$ on the Hilbert space~$\cH$ defined as follows: 
\begin{itemize}
\item $\cH:=\{f\colon G\to\cH_0\mid f(gh)=V_h^*f(g) \text{ if }g\in G,\ h\in G_0\}$; 
\item $U\colon G\to B(\cH)$, $(U_gf)(\tilde{g}):=f(g^{-1}\tilde{g})$; 
\item $\pi\colon\cA\to B(\cH)$, $(\pi(x)f)(g):=\pi_0(\alpha_g^{-1}(x))f(g)$ 
(multiplication operators).
\end{itemize}
\end{defn}

It is useful to reformulate the above definition as follows, in the spirit of \cite[Sect. 3.3]{Se77}. 

\begin{rem}\mlabel{defind}
Let $U\colon G\to B(\cH)$ be a unitary representation of a finite group. 
Assume that a subgroup $G_0\subseteq G$ and a closed linear subspace $\cH_0\subseteq\cH$ 
satisfy $U_{G_0}\cH_0\subseteq\cH_0$, and define $V\colon G_0\to B(\cH_0)$, $V_h:=U_h\vert_{\cH_0}$. 

The group representation $U$ is unitarily equivalent to the representation 
induced from $V$   
if for some (and hence for every) system $\1=g_0,g_1,\dots,g_r$ 
of representatives of the left $G_0$-cosets in $G$ we have an 
orthogonal direct sum decomposition 
\begin{equation}\label{defind_eq1}
\cH=\bigoplus_{j=0}^r U_{g_j}\cH_0. 
\end{equation}
If, moreover, $(\cA,G,\alpha)$ is a $C^*$-dynamical system 
and one has a covariant representation $(\pi_0,V)$ of $(\cA,G_0,\alpha_0)$ on $\cH_0$, 
where $\alpha_0:=\alpha\vert_{G_0}$ and $V\colon G_0\to B(\cH_0)$ is as above, 
then the corresponding induced covariant representation is unitarily equivalent to 
the induced covariant representation $(\pi,U)$ of $(\cA,G,\alpha)$ on $\cH$  
with $\pi\colon\cA\to B(\cH)$ given by 
\[ \pi(a)\vert_{U_{g_j}\cH_0}:=U_{g_j}\pi_0(a)U_{g_j}^{-1}\vert_{U_{g_j}\cH_0}
\quad \mbox{ for } \quad j=0,1,\dots,r, \ x\in\cA.\] 
\end{rem}

\begin{rem}\mlabel{fixes}
For later use, we note that in Definition~\ref{Def3.1} there exists a natural unitary operator 
between the spaces of fixed points 
$$\cH_0^{G_0}\to\cH^G,\quad v\mapsto f_v $$
where $f_v\colon G\to \cH_0$, $f_v(g):=\vert G/G_0\vert^{-1/2}v$ for all $g\in G$ and $v\in \cH_0^{G_0}$. 
\end{rem}

\begin{rem}\mlabel{maxi}
In the setting of Remark~\ref{defind}, let $p_0 \: \cH \to \cH$ denote the 
orthogonal projection onto $\cH_0$. Then 
\[ G_0
=\{g\in G\mid U_g p_0 U_g^{-1}=p_0\} \subseteq\{g\in G\mid \pi_0\circ\alpha_g\simeq\pi_0\} 
\subseteq\{g\in G\mid \pi_0\circ\alpha_g\approx\pi_0\} \]
recalling that $\approx$ stands for quasi-equivalence of $*$-representations. 
\end{rem}

\subsection{Commutants of induced representations}

The situation when both inclusions in Remark~\ref{maxi} are equalities should be thought of 
as a maximality property of the subgroup~$G_0$ in the same way as the polarizations at points 
of duals of finite dimensional Lie algebras are maximal isotropic subalgebras, 
which leads to irreducibility of the corresponding induced representations.  
That idea will be made precise in Corollary~\ref{TakTh4.3cor} below. 
To that end we need 
the following variant of \cite[Th.~4.3]{Tak67} for finite groups acting on $C^*$-algebras that need not be separable. 

\begin{prop}\mlabel{TakTh4.3} 
Let $(\pi,U)$ be the covariant representation of $(\cA,G,\alpha)$ on $\cH$ 
induced by the covariant representation  $(\pi_0,V)$ of $(\cA,G_0,\alpha\vert_{G_0})$ on $\cH_0$.  
Pick any complete system $\1=g_0,g_1,\dots,g_r$ 
of representatives of the left cosets of $G_0$ in $G$
and for $j=0,\dots,r$ 
let $p_j\in B(\cH)$ be the orthogonal projection 
on $\cH_j:= U_{g_j}\cH_0$ and $G_j:=g_jG_0g_j^{-1}$.  

Then $G_j=\{g\in G\mid U_gp_jU_g^{-1}=p_j\}$, 
and if we define $(\pi_j,V_j)$ as the covariant representation of $(\cA,G_j,\alpha\vert_{G_j})$ 
on $\cH_j$, given by $\pi_j(\cdot):=\pi(\cdot)\vert_{\cH_j}$ and $V_j(\cdot):=U(\cdot)\vert_{\cH_j}$, 
then 
$(\pi,U)$ is induced by $(\pi_j,V_j)$ and the map 
$$\Phi_j\colon (\pi{\rtimes U})(\cA\rtimes G)'\cap\{p_0,p_1,\dots,p_r\}'\to (\pi_j{\rtimes V_j})(\cA\rtimes G_j)', 
\quad T\mapsto T\vert_{\cH_j}$$ 
is well defined and is a $*$-isomorphism of von Neumann algebras.
\end{prop}

\begin{prf}
The action of $G$ permutes the subspaces $\cH_0, \ldots,\cH_r$, 
and this implies that  
$(\pi,U)$ is induced by $(\pi_j,V_j)$ 
for $j=0,\dots,r$. 

For every $T\in (\pi{\rtimes U})(\cA\rtimes G)'\cap\{p_0,p_1,\dots,p_r\}'$ we have $Tp_j=p_jT$. 
In particular, such an operator $T$ is given by a block diagonal matrix 
with respect to the decomposition $p_0+\cdots+p_r=\1$, 
and this also shows that for any $j\in\{0,\dots,r\}$ the map $\Phi_j$ 
is a $*$-homomorphism. 

To check that $\Phi_j$ is injective, recall 
that the group $G$ acts transitively on $\{\cH_0,\dots,\cH_r\}$. 
Therefore the subspace $\cH_j$ of $\cH$ is $G$-cyclic, hence separates 
the commutant $U_G'$ which contains 
 $(\pi \rtimes U)(\cA \rtimes G)'$. 

To check that the image of $\Phi_j$ is $(\pi_j{\rtimes V_j})(\cA\rtimes G_j)'$, 
let $T_j\in (\pi_j{\rtimes V_j})(\cA\rtimes G_j)'$. 
For any $g \in G$ with $U_g\cH_j = \cH_k$, we then obtain an operator 
\[ T_k := U_g T_j U_g^* \: \cH_k \to \cH_k \] 
which does not depend on the choice of $g$ in the $G_j$-coset
$\{ x \in G \: U_x\cH_j = \cH_k\}$. 
We thus obtain a block-diagonal operator $T\in B(\cH)$ 
with $T\vert_{\cH_k}:=T_k$ for $k=0,\dots,r$. 
It is straightforward to verify that $T$ 
commutes with $(\pi \rtimes U)(\cA \rtimes G)$. 
Clearly also $T\in\{p_0,p_1,\dots,p_r\}'$ and $\Phi_j(T)=T_j$. 
\end{prf}

\begin{cor}\mlabel{TakTh4.3cor}
Suppose that, in the setting of {\rm Proposition~\ref{TakTh4.3}}, 
the following additional conditions are satisfied: 
\begin{itemize}
\item $G_0=\{g\in G\mid \pi_0\circ\alpha_g\approx\pi_0\}$; 
\item $\pi_0\colon \cA\to B(\cH_0)$ is a factor representation.  
\end{itemize}
Then the map 
$$\Phi_j\colon (\pi \rtimes U)(\cA\rtimes G)'\to 
(\pi_j\rtimes U_j')(\cA\rtimes G_j)', 
\quad T\mapsto p_j T\vert_{\cH_j}$$ 
is a well-defined $*$-isomorphism of von Neumann algebras for $j=0,\dots,r$.
\end{cor}

\begin{prf}
Using Proposition~\ref{TakTh4.3}, it suffices to show that 
$(\pi{\rtimes U})(\cA\rtimes G)'\subseteq\{p_0,p_1,\dots,p_r\}'$. 
One has $(\pi{\rtimes U})(\cA\rtimes G)'=(\pi(\cA)\cup\pi(G))'=\pi(\cA)'\cap\pi(G)'\subseteq\pi(\cA)'$ 
and we will check that for arbitrary $T\in\pi(\cA)'$ one has $Tp_j=p_jT$ for $j=0,\dots,r$. 

In fact, writing the operators on $\cH$ as block matrices with respect to the decomposition 
$p_0+\cdots+p_r=\1$, one has $T=(T_{jk})_{0\le j,k\le r}$ 
and 
$\pi(x)=\diag(\pi_0(x),\pi_1(x),\dots,\pi_r(x))$ is a block diagonal matrix for all $x\in\cA$. 
Hence the condition $T\in\pi(\cA)'$ is equivalent to the fact that for all $j,k=0,\dots,r$ one has 
$$(\forall x\in\cA)\quad \pi_j(x)T_{jk}=T_{jk}\pi_k(x).$$
On the other hand the hypotheses imply that if $j\ne k$ then $\pi_j\vert_{\cA}$ and $\pi_k\vert_{\cA}$ are factor representations that are not quasi-equivalent, hence it follows by \cite[Cor. 5.3.6]{Dix64} that there exists no non-zero intertwiner 
between them, and then $T_{jk}=0$. 
Therefore $T$ is given by a block-diagonal matrix, that is, $Tp_j=p_jT$ for $j=0,\dots,r$, 
and this completes the proof. 
\end{prf}

\begin{rem}\mlabel{well}
In the setting of Corollary~\ref{TakTh4.3cor}, 
it follows that the induced covariant representation~$(\pi,U)$ is irreducible 
or is a factor representation of some type I, II, or III, 
if and only if so is the input representation~$(\pi_0,V)$. 
\end{rem}

\begin{rem}
In the proof of Proposition~\ref{TakTh4.3}, the idea that  
$T=\Phi_0(T)+\cdots+\Phi_r(T)$ 
is a sum of mutually unitarily equivalent operators is related to \cite[Lemma 3]{Co75}. 

On the other hand, since the action of the group $G$ on the set 
$\{p_0,\dots,p_r\}$ is transitive,  
and the isotropy groups $G_0,\dots,G_r$ are mapped to each other by suitable inner automorphisms of $G$, 
this agrees with the uniqueness assertion on $G_0$ in \cite[Th. 5.2]{Tak67}. 
\end{rem}

\subsection{Covariant factor 
representations as induced representations}

We now turn to a question that is somehow converse to what we did so far, 
namely the question of when a given covariant representation is an induced 
representation (cf.\ Proposition~\ref{fact}.)  
We will need the following observation which is essentially a 
byproduct of the proof of \cite[Th. 3.1]{AL10}. 
We recall that the action of $G$ on $\cA$ is called \emph{ergodic} if 
$\cA^G =\C\1$.

\begin{lem}\mlabel{ergo}
If the action of the finite group 
$G$ on $\cA$ is ergodic, then $\dim\cA\le\vert G\vert<\infty$. 
\end{lem}

\begin{prf}
Define 
$$\alpha\colon\ell^1(G)\to B(\cA), \quad 
\alpha(\phi)=\frac{1}{\vert G\vert}\sum_{g\in G}\phi(g)\alpha_g.$$
Let $\pi_1,\dots,\pi_n$ be a system of representatives 
of all equivalence classes of irreducible representations of~$G$. 
For $j=1,\dots,n$ denote by $d_j$ the dimension of $\pi_j$,  
by $\chi_j\colon G\to\C$ the character of $\pi_j$, and define $\cA_j:=\alpha(\chi_j)\cA\subseteq\cA$, 
so that 
$\cA=\cA_1\oplus\cdots\oplus\cA_n.$
On the other hand, since $G$ is in particular a compact group and its action on~$\cA$ is ergodic 
by hypothesis, 
it follows by \cite[Prop. 2.1]{HLS81} 
that $\dim\cA_j\le d_j^2$ for $j=1,\dots,n$. 
Therefore $\dim\cA\le\sum\limits_{j=1}^n d_j^2=\vert G\vert$, 
by Burnside's Theorem, and this concludes the proof. 
\end{prf}

\begin{lem}\mlabel{finite}
Let $(\pi,U)$ be a covariant representation of $(\cA,G,\alpha)$ on $\cH$  
and consider the action of $G$ on $B(\cH)$ 
by $(g,T)\mapsto U_g T U_g^{-1}$. 
\begin{enumerate}[\rm(i)]
 \item\label{finite_item1} If $(\pi,U)$ is a factor representation, then $G$ acts 
ergodically on the abelian von Neumann algebra 
$\cZ:=\pi(\cA)''\cap \pi(\cA)'$  
and $\dim\cZ\le\vert G\vert<\infty$.
\item\label{finite_item2} If $(\pi,U)$ is irreducible, then $G$ acts ergodically on 
$\pi(\cA)'$  
and $\dim\pi(\cA)'\le\vert G\vert<\infty$.
\end{enumerate}
\end{lem}

\begin{prf}
If $(\pi,U)$ is a factor representation, then one has 
$$\begin{aligned}
\cZ^G
&=\pi(\cA)''\cap\pi(\cA)'\cap U_G'=\pi(\cA)''\cap(\pi(\cA)\cup U_G)' \\
&\subseteq(\pi(\cA)\cup U_G)''\cap(\pi(\cA)\cup U_G)' 
=\C\1
\end{aligned}$$
where the latter equality follows by the hypothesis that $(\pi,U)$ is a factor representation. 
Thus $\cZ^G=\C\1$, and then $\dim\cZ\le\vert G\vert$ by Lemma~\ref{ergo}.

If the representation $\pi$ is irreducible, then $\pi(\cA)'\cap U_G'=\C\1$, 
hence the action of $G$ on $\pi(\cA)'$ is ergodic, and then $\dim\pi(\cA)'\le\vert G\vert$ 
by Lemma~\ref{ergo} again.
\end{prf}

We now give a basic result on $*$-representations with finite dimensional commutant, 
for later use along with Lemma~\ref{finite}. 

\begin{lem}\mlabel{commut}
Let $S$ be any $*$-semigroup with a $*$-representation $\pi\colon S\to B(\cH)$.    
If the center of $\pi(S)'$ is finite dimensional, 
then there exist a finite orthogonal direct sum decomposition 
$\cH=\bigoplus\limits_{j=1}^n\cH_j$
with the following properties:
\begin{enumerate}[\rm(i)]
\item\label{commut_item1} 
The orthogonal projections of $\cH$ onto $\cH_1,\dots,\cH_n$, respectively, 
are the minimal central projections  
in~$\pi(S)'$. 
\item\label{commut_item2} 
The representations $\pi_j:=\pi(\cdot)|_{\cH_j}\colon S\to B(\cH_j)$ 
for 
$j=1,\dots,n$, are pairwise non-quasi-equivalent factor representations. 
\item\label{commut_item3} 
If $\pi(S)'$ is also finite dimensional,  
then, for $j=1,\dots,n$, there exist mutually inequivalent 
irreducible $*$-representations $(\rho_j, V_j)$ of $S$ and $k_j \in \N_0$ such that 
$\pi_j \cong \rho_j^{\oplus k_j}$. 
\end{enumerate}
\end{lem}

\begin{prf}
Let $\cM:=\pi(S)'$ and $p_1,\dots,p_n\in\cM$ be the minimal central projections. 
We show that the assertions hold with $\cH_j:=p_j(\cH)$ for $j=1,\dots,n$. 

For Assertion~\eqref{commut_item2}, let $i,j\in\{1,2,\dots,n\}$ with $i\ne j$. 
Since $\pi_i$ and $\pi_j$ are factor representations, 
it follows just as in \cite[Prop. 5.2.9]{Dix64} that it suffices to prove 
that $\pi_i$ and $\pi_j$ re disjoint $*$-representations, 
that is, their only intertwiner is zero. 
Since the set of intertwiners corresponds to 
$p_i \cM p_j$, this follows from $p_i p_j = 0$. 

For Assertion~\eqref{commut_item3}, if $\cM$ is finite dimensional von Neumann algebra, 
then the ideals $p_j \cM$ are isomorphic to matrix algebras $M_{k_j}(\C)$. 
Accordingly, the factor representation $\pi_j$ is type~I and a 
$k_j$-fold multiple of an irreducible representation $(\rho_j,V_j)$. 
\end{prf}

The following result is essentially a special case of \cite[Th. 5.1--5.2]{Tak67}. 
We note that the finite dimensionality of $\cZ$ follows from 
Lemma~\ref{finite}(i). 

\begin{lem}\mlabel{fact}
Let $(\pi,U)$ be a covariant factor representation of $(\cA,G,\alpha)$ on $\cH$ 
and $\cZ:=\pi(\cA)''\cap \pi(\cA)'$.  
Then the following assertions hold: 
\begin{enumerate}[\rm(i)]
\item\label{fact_item2} 
Let $p \in \cZ$ be a minimal projection and $\cK := p\cH$, 
$H:=\{g\in G\mid U_g p U_g^{-1}=p\}$,  
and let $(\rho,V)$ be the covariant representation of $(\cA,H,\alpha\vert_{H})$ 
on $\cK$ given by 
$\rho(a) := \pi(a)\vert_\cK$ and $V_h:=U_h\vert_{\cK}$. 
Then 
$H=\{g\in G\mid \rho\circ\alpha_g\approx\rho\}$, 
$\rho$ is a factor $*$-representation and 
$(\pi,U)$ is unitarily equivalent to the covariant representation induced by 
$(\rho,V)$. 
\item\label{fact_item3} 
The map 
$\Phi\colon (\pi{\rtimes U})(\cA\rtimes G)'\to (\rho{\rtimes V})(\cA\rtimes H)'$, 
$T\mapsto pT\vert_{\cK}$, 
is well defined and is a $*$-isomorphism of von Neumann algebras. 
\end{enumerate}
\end{lem}

\begin{prf} For \eqref{fact_item2}, we note that 
Lemma~\ref{commut}\eqref{commut_item2} implies that $\rho$ 
is a factor representation  
since $\cZ$ is the center of $\pi(\cA)'$. 
Moreover, for every $g\in G$, one has 
$$(\pi\circ\alpha_g)(\cdot)=U_g\pi(\cdot)U_g^{-1},$$ 
so the operators $U_g$ act by conjugation on the von Neumann algebra 
$\pi(\cA)''$, hence also on its center $\cZ$. 
The action of $G$ on $\cZ$ is ergodic by Lemma~\ref{finite}\eqref{finite_item1}, 
and this implies that the action of~$G$ on the set $\{p_1,\dots,p_n\}$ 
of minimal idempotents in $\cZ$ is transitive (cf.~\cite[Lemma 3.4]{Ka13b}). 

Accordingly, $G$ permutes the subspaces $\cH_j := p_j\cH$ transitively, 
and this implies by Remark~\ref{defind} that  
$(\pi,U)$ is unitarily equivalent to the covariant representation 
induced by $(\rho,V)$. 
Since the representations of $\cA$ on the subspaces $\cH_j$ 
are pairwise non-quasi-equivalent by Lemma~\ref{commut}\eqref{commut_item2}
and the group $G$ acts transitively on $\{\cH_1,\dots,\cH_n\}$, 
we obtain 
\[ H=\{g\in G\mid \rho\circ\alpha_g\approx\rho\}.\]
For \eqref{fact_item3}, note that 
$(\pi{\rtimes U})(\cA\rtimes G)''\supseteq\pi(\cA)''\supseteq\cZ$, hence 
using the von Neumann Bicommutant Theorem, 
$(\pi{\rtimes U})(\cA\rtimes G)'=(\pi{\rtimes U})(\cA\rtimes G)'''\subseteq\cZ'=\{p_1,\dots,p_n\}'$, 
and now the conclusion follows by Corollary~\ref{TakTh4.3cor}. 
\end{prf}

\begin{defn}\mlabel{classes}
For any $C^*$-algebra~$\cY$, we denote by $\Fact(\cY)$ the class of 
all unitary equivalence classes $[\pi]$ of 
its factor representations~$\pi$.   

For the $C^*$-dynamical system $(\cA,G,\alpha)$, one has a natural action 
$$G\times \Fact(\cA)\to \Fact(\cA),\quad (g,[\pi])\mapsto [\pi\circ\alpha_{g^{-1}}].$$ 
We define the class $\cC(\cA,G,\alpha)$ of unitary equivalence classes $[(\pi_0,V)]$
of covariant  
representations 
$(\pi_0,V)$ of $C^*$-dynamical systems of the form $(\cA,G_0,\alpha\vert_{G_0})$, 
where $[\pi_0]\in\Fact(\cA)$, $[\pi_0\rtimes V]\in\Fact(\cA\rtimes G_0)$ and 
$G_0$ is the isotropy group in $G$ of the quasi-equivalence class of~$\pi_0$. 
There is a natural action 
\begin{equation}\label{quasi_eq1}
G\times \cC(\cA,G,\alpha)\to\cC(\cA,G,\alpha),\quad 
(g,[(\pi_0,V)])\mapsto [(\pi_0\circ\alpha_{g^{-1}},V\circ c_g^{-1}\vert_{gG_0g^{-1}})]
\end{equation}
where $c_g\colon G\to G, h\mapsto ghg^{-1}$. 
This makes sense since $gG_0g^{-1}$ is the isotropy group of the quasi-equivalence class 
of the factor representation $\pi_0\circ\alpha_{g^{-1}}$ of $\cA$, 
and the covariant representation $(\pi_0\circ\alpha_{g^{-1}},V\circ c_g^{-1}\vert_{gG_0g^{-1}})$ 
is a factor representation. 

It is easily checked that 
two unitarily equivalent covariant representations 
give rise by induction to unitarily equivalent representations, 
and it then follows by Corollary~\ref{TakTh4.3cor} and Remark~\ref{well} 
that the induction of covariant representations defines a map 
\begin{equation}\label{quasi_eq2}
\Theta\colon \cC(\cA,G,\alpha)\to \Fact(\cA\rtimes_\alpha G).
\end{equation}
\end{defn}

\begin{thm}\mlabel{factbij}
For any $C^*$-dynamical system $(\cA,G,\alpha)$ for which the group $G$ is finite, 
the map $\Theta$ of \eqref{quasi_eq2} 
is constant on the orbits of the group action~\eqref{quasi_eq1} 
and defines a bijective correspondence 
\[ \oline\Theta \: \cC(\cA,G,\alpha)/G\to  \Fact(\cA\rtimes_\alpha G)\] 
which preserves $*$-isomorphism classes of commutants, irreducibility, 
and the types I,II,III of factor representations. 
\end{thm}

\begin{prf}
It follows by Lemma~\ref{fact} that $\Theta$ is surjective. 

Injectivity of $\oline\Theta$ follows by Lemma~\ref{TakTh4.3}, 
where the subgroup $G_0$ is determined as the centralizer of a minimal 
central idempotent in $\pi(\cA)'$, and this determines the corresponding 
factor representation $\pi_0$ of $\cA$. 
By Remark~\ref{well} the map $\Theta$ preserves $*$-isomorphism classes of commutants, and we are done. 
\end{prf}

\begin{rem}\mlabel{nonconj}
Resume the setting of Definition~\ref{classes} 
and denote by $\cG(\cA,G,\alpha)$ the set of all subgroups $G_0\subseteq G$ for which 
there exists $[\pi_0]\in\Fact(\cA)$ such that $G_0$ is the isotropy group of the equivalence class of $\pi_0$. 
If we consider the natural action of $G$ on $\cG(\cA,G,\alpha)$ by conjugation, 
then one obtains a $G$-equivariant correspondence 
$\cC(\cA,G,\alpha)\to\cG(\cA,G,\alpha)$ 
that associates to every 
$[(\pi_0,V)]$ in $\cC(\cA,G,\alpha)$ the domain of definition of $V$. 
Using Theorem~\ref{factbij}, this leads to a surjective correspondence  
$$\Fact(\cA\rtimes_\alpha G)\simeq\cC(\cA,G,\alpha)/G\to\cG(\cA,G,\alpha)/G$$
which in particular defines an equivalence relation on $\Fact(\cA\rtimes_\alpha G)$ 
with finitely many equivalence classes. 

This remark provides the following necessary criterion 
for unitary equivalence of two covariant factor representations $(\pi_1,U_1)$ and $(\pi_2,U_2)$ of $(\cA,G,\alpha)$: 
For $j=1,2$, assume that $\pi_j$ is induced from some factor covariant representation $(\pi_{j0},V_j)$ of 
$(\cA,G_j,\alpha\vert_{G_j})$ with  
$[\pi_{j0}]\in\Fact(\cA)$ and 
$G_j$ the isotropy group of the quasi-equivalence class of $\pi_{j0}$. 
If the subgroups $G_1$ and $G_2$ fail to be conjugated to each other in $G$, 
then the covariant factor representations $(\pi_1,U_1)$ and $(\pi_2,U_2)$ are not 
unitarily equivalent, 
and equivalently, $[\pi_1\rtimes U_1]\ne[\pi_2\rtimes U_2]$ in $\Fact(\cA\rtimes G)$. 
\end{rem}

\begin{ex}\mlabel{regular}
In the setting of Definition~\ref{classes}, consider the subclass $\cC_0(\cA,G,\alpha)$ 
that corresponds to $G_0=\{\1\}$. 
Hence the elements of $\cC_0(\cA,G,\alpha)$ can be identified with the unitary equivalence classes 
of factor $*$-representations $\pi_0\colon\cA\to B(\cH_0)$ 
with the property that for every $\1 \not=g\in G$ the factor representation 
$\pi_0\circ\alpha_g$ is not quasi-equivalent to $\pi_0$. 
Then the corresponding induced representation $(\pi,U)$ of $(\cA,G,\alpha)$ 
acts on the Hilbert space $\cH = \ell^2(G,\cH_0)$ 
consisting of all functions $f\colon G\to\cH_0$ 
and is given by the formulas from Definition~\ref{Def3.1}. 
The image of the equivalence class $[\pi_0]$ in $\cC(\cA,G,\alpha)/G$ 
is ${\{\pi_0\circ\alpha_g\mid g\in G\}}$, and all the representations from this set induce 
the same element of $\Fact(\cA\rtimes_\alpha G)$ by Theorem~\ref{factbij}. 

In particular, this shows that there is no analogue of Theorem~\ref{factbij} with unitary equivalence classes 
of factor representations replaced by quasi-equivalence classes. 
\end{ex}

\begin{ex}\mlabel{irregular}
To illustrate Definition~\ref{classes} by a case which is opposite to Example~\ref{regular}, 
consider the subclass $\cC_1(\cA,G,\alpha)$ 
that corresponds to $G_0=G$. 
Hence the elements of $\cC_1(\cA,G,\alpha)$ are the unitary equivalence classes 
of covariant factor representations $(\pi_0,V)$, 
where $\pi_0\colon\cA\to B(\cH_0)$ is also factor $*$-representation. 
As $(\pi_0,V)$ is a covariant representation of $(\cA,G,\alpha)$, 
it follows that, for every $g\in G$, the factor representation 
$\pi_0\circ\alpha_g$ is quasi-equivalent (actually unitarily equivalent) to $\pi_0$. 
Then the corresponding induced representation $(\pi,U)$ of $(\cA,G,\alpha)$ 
is given by the formulas from Definition~\ref{Def3.1} 
and is easily seen to be unitarily equivalent to $(\pi_0,V)$. 
The image of the equivalence class $[(\pi_0,V)]$ in $\cC(\cA,G,\alpha)/G$ 
is the singleton $\{[(\pi_0,V)]\}$, and this corresponds to $[\pi_0\rtimes V]\in \Fact(\cA\rtimes_\alpha G)$ via  Theorem~\ref{factbij}. 

For later use we also note a simple necessary criterion for unitary equivalence 
of factor representations $[\pi_{0j}\rtimes V_j]$, j=1,2,  
that belong to the image of $\cC_1(\cA,G,\alpha)$ in $\Fact(\cA\rtimes_\alpha G)$: 
if $[\pi_{01}\rtimes V_1]=[\pi_{02}\rtimes V_2]$, then $[\pi_{01}]=[\pi_{02}]\in\Fact(\cA)$ 
and the unitary representations $V_1$ and $V_2$ of $G$ are unitarily equivalent. 
\end{ex}

\subsection{Classification of irreducible covariant representations} 

In the setting of Theorem~\ref{factbij}, the definition of the class $\cC(\cA,G,\alpha)$ 
is somewhat involved, in the sense that 
for a covariant representation $(\pi_0,V)$ of a $C^*$-dynamical system $(\cA,G_0,\alpha\vert_{G_0})$, 
the unitary representation $V$ of $G_0$ is rather implicitly described 
by the set of conditions $\pi_0\in\Fact(\cA)$, $\pi_0\rtimes V\in\Fact(\cA\rtimes G_0)$,  
$G_0$ is the isotropy group of the quasi-equivalence class of $\pi_0$. 
We will now discuss this issue for the subclass of irreducible covariant representation $(\pi_0,V)$, 
thus recovering the precise description of the dual space of $\cA\rtimes G$ 
as in \cite{Tak67}, \cite{AL10}, and \cite{Ka13}. 

We will basically classify the factor representation $\rho$ from Lemma~\ref{fact} 
in the case when the covariant representation~$(\pi,U)$ is irreducible. 
See also \cite[Th. 3.1]{AL10} and the comments after \cite[Th. 4.2]{Ka13b}. 
If $S$ is any group, $\cH$ is any Hilbert space, and $V\colon S\to B(\cH)$ 
and $\sigma\colon S\times S\to\T$ are maps satisfying 
the $V_{ts}=\sigma(t,s)V_t V_s$ for all $s,t\in S$, 
then we say that $V$ is a \emph{$\sigma$-projective representation}.  

For the following proposition, we recall the existence of the minimal 
central projections from Lemmas~\ref{finite} and \ref{commut}.

\begin{prop}\mlabel{irr}
Let $(\pi,U)$ be an irreducible covariant representation of $(\cA,G,\alpha)$ on $\cH$. 
Let $p_1,\dots,p_n$ be the minimal central projections in $\pi(\cA)'$. 
Then the following assertions hold: 
\begin{enumerate}[\rm(i)]
\item\label{irr_item1} 
For $j\in\{1,\dots,n\}$ let $d_j^2:=\dim (\pi(\cA)'p_j)$ 
and $G_j:=\{g\in G\mid U_g p_j U_g^{-1}=p_j\}$.   
Then  $\pi_j\cong \rho_j^{\oplus d_j}$ for an 
irreducible representation $(\rho_j, \cK_j)$ of $\cA$ 
and the covariant representation $(\pi_j,V_j)$ of $(\cA,G_j,\alpha\vert_{G_j})$ 
on $\cH_j$ is irreducible. 
\item\label{irr_item2} 
$G_j=\{g\in G\mid \rho_j\circ\alpha_g\simeq\rho_j\}$ and $d_j=d_1$ for $j=1,\dots,n$. 
\item\label{irr_item3} 
The covariant representation $(\pi_1,V_1)$ is unitarily equivalent to a 
covariant representation 
$(\rho_1\otimes\1,V_1^{(1)} \otimes V_1^{(2)})$, 
where $V_{1}^{(1)}\colon G_1\to \U(\cK_1)$ 
is a unitary $\sigma$-projective representation 
for a $2$-cocycle $\sigma\colon G_1\times G_1\to\T$ whose cohomology 
class is uniquely determined, and 
$V_{1}^{(2)}\colon G_1\to \U_{d_1}(\C)$ 
is a unitary irreducible $\overline{\sigma}$-projective representation. 
\item\label{irr_item4} 
The map $\Theta$ of \eqref{quasi_eq2} defines by restriction a bijection 
onto $(\cA\rtimes G)\,\hat{}\ $ 
from the $G$-orbits of equivalence classes 
 $[(\rho_{1}\otimes\1,V_{1}^{(1)}\otimes V_{1}^{(2)})]$ 
as above, where $V_{1}^{(2)}$ can be any irreducible $\overline{\sigma}$-projective unitary representation 
of the isotropy group $G_1$ of $[\rho_{1}]\in\widehat{\cA}$, 
and unitary inequivalent choices of $V_{1}^{(2)}$ lead to distinct equivalence classes 
 $[(\rho_{1}\otimes\1,V_{1}^{(1)}\otimes V_{1}^{(2)})]$. 
\end{enumerate}
\end{prop}

\begin{prf} (i) The covariant representation $(\pi_j,V_j)$ 
of $(\cA,G_j,\alpha\vert_{G_j})$ on $\cH_j$ is irreducible 
by Proposition~\ref{fact}\eqref{fact_item3}.  

(ii) follows by Lemma~\ref{commut}\eqref{commut_item3}. 

(iii) follows by \cite[Lemmas 5.1-5.2]{Tak67}, 
where the $C^*$-algebra needs not be separable if the group is finite.  
In fact, in this case, \cite[Th. 2.6]{Tak67} (needed in the proof of \cite[Lemma 5.1]{Tak67}) 
is actually the trivial observation that, 
because 
\[ G_1=\{g\in G\mid \rho_{1}\circ\alpha_g\simeq\rho_{1}\},\]
there exist a $2$-cocycle 
$\sigma\colon G_1\times G_1\to\T$, 
and a unitary $\sigma$-projective representation \break $V_{1}^{(1)}\colon G_1\to 
\U(\cK_{1})$,  
with 
\[ (\rho_{1}\circ\alpha_g)(x)=V_{1}^{(1)}(g)\rho_{1}(x)V_{1}^{(1)}(g)^* 
\quad \mbox{ for } \quad g\in G_1, x\in\cA.\] 
The cohomology class $[\sigma]$ is uniquely determined because the 
corresponding central $\T$-extension of $G_1$ is isomorphic to the 
pullback of the $\T$-extension of $\PU(\cK_1):=\U(\cK_1)/\T\1$ under the corresponding 
projective representation. 
Then one has 
\[ (V_{1}^{(1)}(g)\otimes\1)\pi_1(x)(V_{1}^{(1)}(g)\otimes\1)^*=(\pi_1\circ\alpha_g)(x)=V_g\pi_1(x)V_g, \] 
hence $(V_{1}^{(1)}(g)\otimes\1)^*V_g$ belongs to the commutant of 
$\pi_1\cong \rho_1^{\oplus d_1}$, 
hence there exists a $V_{1}^{(2)}\in \U_{d_1}(\C)$ with $(V_{1}^{(1)}(g)\otimes\1)^*V_g=\1\otimes (V_{1}^{(2)}(g)$. 
Since $V$ is a unitary representation, $V\simeq V_{1}^{(1)}\otimes V_{1}^{(2)}$, and $V_{1}^{(1)}$ is a $\sigma$-projective representation, it follows that $V_{1}^{(2)}$ is a $\overline{\sigma}$-projective unitary representation of $G_1$, 
which is also irreducible, exactly as in the proof of \cite[Lemma 5.2]{Tak67}. 

Conversely, given $[\rho_{1}]\in\widehat{\cA}$, one determines its isotropy group $G_1$ 
with respect to the action of $G$ on $\widehat{\cA}$, then a cocycle 
$\sigma\in Z^2(G_1,\T)$ and a $\sigma$-projective representation $V_{1}^{(1)}\colon G_1\to B(\cK_{1})$, 
and then every irreducible $\overline{\sigma}$-projective unitary representation $V_1^{(2)}$ 
gives rise to a covariant representation $(\rho_{1}\otimes\1,V_{1}^{(1)}\otimes V_{1}^{(2)})$ of $(\cA,G_1,\alpha\vert_{G_1})$ 
which is an irreducible representation and induces an irreducible covariant representation of $(\cA,G,\alpha)$. 
Thus (iv) follows by Theorem~\ref{factbij}, and we are done. 
\end{prf}

\subsection{Representations of fixed-point subalgebras}

We will need the following slight generalization of the well-known fact that factor $*$-re\-pre\-sen\-tations are nondegenerate. 

\begin{lem}\mlabel{factres}
Let $\cB_0$ be a $*$-subalgebra of some $*$-algebra $\cB$. 
Then, for any factor representation $\pi\colon\cB\to B(\cH)$, 
one has $\cB_0\not\subset\ker\pi$ if and only if the subspace $\spann\,\pi(\cB\cB_0)\cH$ is dense in $\cH$. 
\end{lem}

\begin{prf}
The hypothesis that $\pi$ is a factor representation is equivalent to the fact that there exists 
no nontrivial closed linear subspace of $\cH$  that it is invariant both under 
$\pi(\cB)$ and its commutant $\pi(\cB)'$, 
since the orthogonal projection on such a subspace would belong to the center 
$\pi(\cB)'\cap\pi(\cB)''$ of $\pi(\cB)''$. 
 
We may assume $\cH\ne\{0\}$ and the assertion will follow by the above observation 
as soon as we will have proved that if $\cB_0\not\subset\ker\pi$ then $\pi(\cB\cB_0)\cH\ne\{0\}$, 
since this subset of $\cH$ is invariant under $\pi(\cB)$ and $\pi(\cB)'$. 
In fact, since $\cB_0\not\subset\ker\pi$, there exists $b_0\in\cB_0$ 
with $\pi(b_0)\ne 0$ and then 
$\pi(b_0^*) \pi(b_0)\cH \subeq \pi(\cB)\pi(\cB_0)\cH$ is non-zero. 
\end{prf}

\begin{prop}\mlabel{factim}
Let $\cB_0\subseteq\cB$ be a $*$-subalgebra of some Banach $*$-algebra with approximate unit, 
with $\cB_0\cB\cB_0\subseteq\cB_0$. 
Denote by $\Fact_{\cB_0}(\cB)$ the class of all $[\pi]\in\Fact(\cB)$ with $\cB_0\not\subset\ker\pi$ 
and denote by $\Psi([\pi])$ the unitary equivalence class of the 
$\cB_0$-nondegenerate part of $\pi\vert_{\cB}$ 
for $[\pi]\in \Fact_{\cB_0}(\cB)$. 
Then the correspondence 
$$\Psi\colon\Fact_{\cB_0}(\cB)\to\Fact(\cB_0)$$
is a bijection and a homeomorphism with respect to the regional topologies, 
which preserves $*$-isomorphism classes of commutants, irreducibility, 
and the types I,II,III of factor representations. 
\end{prop}

\begin{prf}
This follows by \cite[vol. II, XI.7.6, Prop.~and Rem.~1]{FD88}, 
using the above Lemma~\ref{factres}. 
\end{prf}

For the following statement we recall that for any $C^*$-algebra $\cB$ we denote by $M(\cB)$ its multiplier algebra, 
which is a unital $C^*$-algebra with a canonical embedding $\cB\hookrightarrow M(\cB)$, 
and one has $\cB=M(\cB)$ if and only if $\cB$ is unital. 

\begin{prop}\mlabel{cross}
Denote $\cA^G:=\{x\in\cA\mid(\forall g\in G)\ \alpha_g(x)=x\}$ 
and define $\widetilde{\cA}:=\cA + \C \1 \subeq M(\cA)$. 
Let  $p\in \ell^1(G,\widetilde{\cA};\alpha)\subseteq M(\cA\rtimes_\alpha G)$ 
be the constant function that is equal to $\1\in \widetilde{\cA}$ at every point of $G$. 
Then the following assertions hold: 
\begin{enumerate}[\rm(i)]
\item\label{cross_item1} 
One has $p=p^*=p^2$, and by restriction and co-restriction the injective linear map  
$\cA\to \cA\rtimes_\alpha G$,  
$x\mapsto\alpha(\cdot,x)$,  
yields a $*$-isomorphism $\iota\colon \cA^G\to p(\cA\rtimes_\alpha G)p$.  
\item\label{cross_item2} 
For every factor covariant representation $(\pi,U)$ of $(\cA,G,\alpha)$,  
we will write $(\pi\rtimes U)\,\tilde{}$ for 
the canonical extension of $\pi{\rtimes U}$ from $\cA\rtimes_\alpha G$ 
to $M(\cA\rtimes_\alpha G)$,  
and by $(\pi\rtimes U)^G$ the representation of $\cA^G$ obtained as 
the nondegenerate part of $(\pi\rtimes U)\,\tilde{}\circ\iota$. 
We also define the following open subset of $\Fact(\cA\rtimes G)$, 
$$
\Fact_p(\cA\rtimes_\alpha G)
:=\{[\pi\rtimes U]\in\Fact(\cA\rtimes_\alpha G)\mid (\pi\rtimes U)\,\tilde{}(p)\ne0\}. $$
Then the map 
$$\Fact_p(\cA\rtimes_\alpha G)\to \Fact(\cA^G),\quad 
[\pi{\rtimes U}]\mapsto[(\pi{\rtimes U})^G]$$ 
is a well-defined homeomorphism which preserves $*$-isomorphism classes of commutants. 
\item\label{cross_item3} 
One has 
$[\pi{\rtimes U}]\in \Fact_p(\cA\rtimes_\alpha G)$
if and only if 
\[ [\pi{\rtimes U}]\in\Fact(\cA\rtimes_\alpha G) \quad \mbox{ and } \quad 
\sum\limits_{s\in G}U_s \ne 0,\] 
and in this case $P_U:=\frac{1}{\vert G\vert}\sum\limits_{s\in G}U_s$ 
is the orthogonal projection onto the essential space of the $*$-representation $(\pi{\rtimes U})\circ\iota$ of $\cA^G$.
\item\label{cross_item4} 
Let $G_0\subseteq G$ be any subgroup and $(\pi_0,V)$ be any covariant representation  
of the $C^*$-dynamical system $(\cA,G_0,\alpha\vert_{G_0})$  
for which the corresponding induced covariant representation 
$(\pi,U)$  
of $(G,\cA,\alpha)$ is irreducible. 
Then $[\pi{\rtimes U}]\in\Fact_p(\cA\rtimes_\alpha G)$ 
if and only if $\sum\limits_{s\in G_0}V_s \ne 0$. 
In particular, this is always the case if $G_0=\{\1\}$.   
\item\label{cross_item5} The $C^*$-algebra $\cA$ is of type I  
if and only if $\cA^G$ is, and if and only if $\cA\rtimes_\alpha G$ is. 
\end{enumerate}
\end{prop}

\begin{prf}
Since $G$ is in particular a compact group, Assertion~\eqref{cross_item1} follows by \cite{Ro79} 
(see also \cite[II.10.4.18]{Bl06}).   

For Assertion~\eqref{cross_item2}, 
denote $\cB:=\cA\rtimes_\alpha G=\ell^1(G,\cA;\alpha)$, 
$\cB_0:=p\cB p$, 
and note that, 
since $p=p^*=p^2\in M(\cB)$, 
$p\cB p$ is a hereditary subalgebra of $\cB$,  
that is, $\cB_0\cB\cB_0\subseteq\cB_0$.  
Then the assertion follows by Proposition~\ref{factim} along with Assertion~\eqref{cross_item1}. 

For Assertion~\eqref{cross_item3}, recall that if $(\pi,U)$ is a covariant representation of $(\cA,G,\alpha)$, 
then the corresponding representation of $\cB=\ell^1(G,\cA;\alpha)$ 
is 
\[ \pi\rtimes U\colon \cB\to B(\cH),\quad 
(\pi\rtimes U)(f)=\frac{1}{\vert G\vert}\sum\limits_{s\in G}\pi(f(s))U_s.\] 
If $(\pi,U)$ is a factor representation, 
then also $\pi\rtimes U\colon \cB\to B(\cH)$ is a factor representation hence it is nondegenerate, 
and then it extends to $M(\cB)\supseteq \ell^1(G,\widetilde{\cA};\alpha)$.  
By the above formula one then obtains for the constant function $p\colon G\to \widetilde{\cA}$, 
$$(\pi\rtimes U)(p)=\frac{1}{\vert G\vert}\sum\limits_{s\in G}U_s.$$
Therefore the condition $(\pi\rtimes U)(p)\ne 0$ is equivalent to $\sum\limits_{s\in G}U_s \ne 0$. 

For Assertion~\eqref{cross_item4}, let $\{\1=g_0,g_1,\dots,g_n\}$ be 
a complete system of representatives of the $G_0$-left cosets sets in~$G$. 
Note that 
\begin{equation*}
\Bigl(\sum\limits_{g\in G}U_g\Bigr)\cH=\{v\in\cH\mid(\forall g\in G)\ U_g v=v\}=:\cH^G
\end{equation*}
and similarly
\begin{equation*}
\Bigl(\sum\limits_{s\in G_0}V_g\Bigr)\cH_0=\{v_0\in\cH_0\mid(\forall s\in G_0)\ V_s v_0=v_0\}=:\cH_0^{G_0}, 
\end{equation*}
hence we must show that $\cH^G\ne\{0\}$ if and only if $\cH_0^{G_0}\ne\{0\}$, 
which follows by Remark~\ref{fixes}. 

Finally, Assertion~\eqref{cross_item4} follows by \cite[Th. 4.1]{Ri80}, and we are done.  
\end{prf} 

Our next aim is to describe the dual space $\widehat{\cA^G}$ using Proposition~\ref{cross}. 
In order to simplify the corresponding statement (Theorem~\ref{fix}) 
we make the following definition. 

\begin{defn}\mlabel{masamichi}
For the $C^*$-dynamical system $(\cA,G,\alpha)$,  
let $\pi_0\colon\cA\to B(\cH_{\pi_0})$ be an arbitrary irreducible $*$-representation 
and 
$$G_{\pi_0}:=\{g\in G\mid \pi_0\circ\alpha_g\simeq \pi_0\}.$$ 
Let $\omega_0\colon G_{\pi_0}\times G_{\pi_0}\to\T$ be a $2$-cocycle  
and $V_0\colon  G_{\pi_0}\to \U(\cH_{\pi_0})$ be any $\omega_0$-projective unitary representation 
 satisfying
$$(\forall g\in G_{\pi_0})(\forall a \in \cA)
\quad V_0(g)\pi_0(a)V_0(g)^*=\pi_0(\alpha_g(a)). $$
Let $\Gamma(\cA,G,\alpha)$ be the class of all 
covariant representations 
$(\pi,V)$ 
of $(\cA,G_{\pi_0},\alpha\vert_{G_{\pi_0}})$ 
obtained as 
$$
\pi:=\pi_0\otimes\id_{\cH_{\overline{U_0}}}
\text{ and } 
V:=V_0\otimes\overline{U_0}$$
where  
$U_0\colon G_{\pi_0}\to B(\cH_{U_0})$ is an arbitrary irreducible $\omega_0$-projective unitary representation, 
and let $\Gamma_0(\cA,G,\alpha)\subseteq\Gamma(\cA,G,\alpha)$ be the subclass 
defined by imposing the additional condition
\begin{equation}\label{fix_eq1}
\sum_{g\in G_{\pi_0}}V_0(g)\otimes\overline{U_0}(g)\ne0. 
\end{equation}
Then the set $\widetilde{\Gamma}_0(\cA,G,\alpha)$ of all unitary equivalence classes in $\Gamma_0(\cA,G,\alpha)$, 
is naturally acted on by $G$ via $\alpha$ and via the conjugation 
 action of $G$ on itself 
as in \eqref{quasi_eq1} (see also \cite[Th.~7.2]{Tak67}, \cite[Sect.~2]{Ka13}), 
and we denote by $\widetilde{\Gamma}_0/G$ the corresponding quotient set. 

Finally, we define 
\begin{equation}
  \label{eq:cI}
\cI\colon\widetilde{\Gamma}_0(\cA,G,\alpha)/G \to\widehat{\cA^G}
\end{equation}
to be the correspondence that takes the unitary equivalence class 
$[(\pi,V)]$ of a covariant representation 
of $(\cA,G_{\pi_0},\alpha\vert_{G_{\pi_0}})$ 
to the unitary equivalence class of $\Ind_{G_{\pi_0}}^G V$
and then composing that correspondence with the map $[\Pi]\mapsto[\Pi^G]$ 
given by Proposition~\ref{cross}\eqref{cross_item2}.
\end{defn}

The following result should be compared with \cite[Th. 3.3]{Ka13} 
(or \cite[Th. 7.2]{Tak67} if the $C^*$-algebra $\cA$ is of type~I). 

\begin{thm}\mlabel{fix}
For every $C^*$-dynamical system $(\cA,G,\alpha)$,  whose group~$G$ is finite,  
the map $\cI\colon\widetilde{\Gamma}_0(\cA,G,\alpha)/G\to\widehat{\cA^G}$ 
is a bijection.  
\end{thm}

\begin{prf}
Using the notation of Definition~\ref{masamichi}, 
one can see that Proposition~\ref{irr} 
 establishes 
the bijection \eqref{eq:cI} from the set of $G$-orbits in $\widetilde{\Gamma}(\cA,G,\alpha)$ onto $(\cA\rtimes G)\,\hat{}\ $.

In some more detail, to check surjectivity of $\cI$, we use our Lemma~\ref{fact} 
to show that every irreducible covariant representation of $(\cA,G,\alpha)$ 
is induced from a factor representation of type~I. 
Now the conclusion follows by Proposition~\ref{cross}\eqref{cross_item2} and \eqref{cross_item4}. 
\end{prf}

\section{Representations of symmetric tensor powers}
\mlabel{sec:5}

It follows by \cite{BN15} that the representations of exponential $C^*$-algebras 
$e^\cA = \oplus_{n \in \N_0}^{c_0} S^n(\cA)$ 
play an important role in the representation theory of ball semigroups. 
Since any exponential $C^*$-algebra is the $c^0$-direct sum of the symmetric 
tensor powers $S^n(\cA) = (\cA^{\otimes n})^{S_n}$ of a $C^*$-algebra, 
one needs to understand the representations of 
the algebras $S^n(\cA)$ in terms of the representations of $\cA$. 

\begin{defn}\mlabel{def11.11}
If $\cA$ is a $C^*$-algebra, then 
for every permutation $\sigma\in S_n$ we denote by $\alpha_\sigma\colon \cA^{\otimes n}\to\cA^{\otimes n}$ 
its corresponding permutation action, and we denote by $(\cA^{\otimes n},S_n,\alpha)$ 
the $C^*$-dynamical system obtained in this way. 
\end{defn}

We record the following fact for later use. 

\begin{lem}\mlabel{lprf0}
If $n\ge 2$, $\cA_1,\dots,\cA_n$ are $C^*$-algebras and  
$\cA:=\cA_1\otimes\cdots\otimes\cA_n$, then the map  
\begin{equation}\label{new_prf0}
\widehat{\cA_1}\times\cdots\times\widehat{\cA_n}\to\widehat{\cA}, \quad 
([\pi_1],\dots,[\pi_n])\mapsto [\pi_1\otimes\cdots\otimes\pi_n]
\end{equation}
is injective, and it is also surjective if at least $n-1$ of the $C^*$-algebras $\cA_1,\dots,\cA_n$ 
are of type~I. 
\end{lem}

\begin{prf}
If $n=2$, the assertion follows by \cite[p.~333]{OT66},  
or \cite[IV.3.4.22/26]{Bl06},
then use induction. 
\end{prf}

\subsection{The case of $C^*$-algebras of type I}

We are now in a position to 
describe the representations of 
the algebras $S^n(\cA)$ in terms of the representations of $\cA$
for any  $C^*$-algebra $\cA$ of type~I, 
by basically specializing Theorem~\ref{fix} to the $C^*$-dynamical system 
$(\cA^{\otimes n},S_n,\alpha)$. 
For the following statement we recall that a \emph{full corner} of a $C^*$-algebra $\cA$ is any subalgebra of the form $p\cA p$, 
where $p=p^2=p^*\in M(\cA)$,  
that is not contained in any closed two-sided ideal of $\cA$, 
(see \cite{Br77}). 

\begin{thm}\mlabel{exact}
Let $\cA$ be a  $C^*$-algebra of type~I and $n\ge 1$. 
For the $C^*$-dynamical system $(\cA^{\otimes n},S_n,\alpha)$,  
define the projection $p\in M(\cA^{\otimes n}\rtimes_\alpha S_n)$ 
and the $*$-isomorphism 
\[ {\iota\colon S^n(\cA)\to p(\cA^{\otimes n}\rtimes_\alpha S_n)p} \] 
as in {\rm Proposition~\ref{cross}\eqref{cross_item1}}. 
Then the following assertions hold: 
\begin{enumerate}[\rm(i)]
\item\label{exact_item1} 
The hereditary subalgebra $p(\cA^{\otimes n}\rtimes_\alpha S_n)p$ is a full corner of $\cA^{\otimes n}\rtimes_\alpha S_n$. 
\item\label{exact_item2} 
For every irreducible representation $\Pi$ of $\cA^{\otimes n}\rtimes_\alpha S_n$ 
one has $\tilde\Pi(p)\ne 0$. 
\item The map $[\Pi]\mapsto[\Pi^G]$ 
from {\rm Proposition~\ref{cross}\eqref{cross_item2}} is a homeomorphism 
from 
$(\cA^{\otimes n}\rtimes_\alpha S_n)\,\hat{}$ onto $S^n(\cA)\,\hat{}\ $. 
\item\label{exact_item3} 
For a partition $q_1+\cdots+q_m=n$, let 
$\bq := (q_1, \ldots, q_m)$ and $S_\bq:=S_{q_1}\times\cdots\times S_{q_m}\subeq S_n$. 
Let $\Gamma_0$ be the family of  
covariant representations $(\Pi_1,W_1)$ of 
$(\cA^{\otimes n},S_{\bq},\alpha\vert_{S_{\bq}})$ 
obtained as 
$$\Pi_1:=\pi_1^{\otimes q_1}\otimes\cdots\otimes\pi_m^{\otimes q_m}
\otimes\id_{\cH_{\overline{U_{01}}}}\otimes\cdots\otimes\id_{\cH_{\overline{U_{0m}}}}
\text{ and } 
W_1:=V_0\otimes\overline{U_{01}}\otimes\cdots\otimes\overline{U_{0m}} $$
where $\pi_k\colon \cA\to B(\cH_{\pi_k})$ 
for $k=1,\dots,m$ are mutually inequivalent irreducible $*$-representations,  
$V_0\colon S_{\bq}\to \U(\cH_{\pi_1}^{\otimes q_1}\otimes\cdots\otimes\cH_{\pi_m}^{\otimes q_m})$ 
is the (external) tensor product of permutation representations, 
 and $U_{0k}\colon S_{q_k}\to \U(\cH_{U_{0k}})$ is an arbitrary irreducible unitary representation 
 for $k=1,\dots,m$. 
Then the set $\widetilde{\Gamma}_0$ of all unitary equivalence classes in $\Gamma_0$, 
is naturally acted upon by $S_n$ via $\alpha$, 
and the bijection 
$\widetilde{\Gamma}_0/S_n\to (\cA^{\otimes n}\rtimes_\alpha S_n)\,\hat{}$  
takes every covariant representation $(\Pi_1,W_1)$  
to its induced covariant representation of $(\cA^{\otimes n},S_n,\alpha)$.  
\end{enumerate}
\end{thm}

\begin{prf} (i) By \cite[Lemma 1]{Br77} this 
is a consequence of (ii). 

(ii) Since $\cA$ is a  $C^*$-algebra of type~I, 
it follows by Lemma~\ref{lprf0} that 
the map \eqref{new_prf0} with ${\cA_1=\cdots=\cA_n=\cA}$ 
is bijective for every $n\ge 1$. 
This implies that, for an irreducible representation $\Pi_0\colon \cA^{\otimes n}\to B(\cH_{\Pi_0})$,  
there exist uniquely determined  
irreducible representations $\pi_1,\dots,\pi_n$ of $\cA$ with $\Pi_0\cong 
\pi_1\otimes\cdots\otimes\pi_n$. After conjugating with a suitable permutation, 
we may assume  
\[ \Pi_0 \cong \rho_1^{\otimes q_1} \otimes \cdots 
\otimes \rho_m^{\otimes q_m},\] 
so that the isotropy group of $[\Pi_0]$ 
is conjugate to $S_{\bq}$. 

Then, in the notation of Theorem~\ref{fix}, one has $\omega_0\equiv 1$ 
and $V_0\colon (S_n)_{\Pi_0}\to \U(\cH_{\Pi_0})$ 
is the (external) tensor product of the permutation representations 
$V_{0k}\colon S_{q_k}\to \U(\cH_{\rho_k}^{\otimes q_k})$. 
Since $\omega_0\equiv 1$, 
it follows that $U_0\colon (S_n)_{\Pi_0}\to \U(\cH_{U_0})$ is an arbitrary irreducible unitary representation 
satisfying~\eqref{fix_eq1}.  
By \cite[Cor. to Th. 2]{OT66},  
there exist uniquely determined irreducible unitary representations 
$U_{0k}\colon S_{q_k}\to  \U(\cH_{U_{0k}})$ for $k=1,\dots,m$ 
whose (external) tensor product is~$U_0$. 
Then, using the fact that, if $T_k\in B(\cH_{0k})$ for $k=1,\dots,m$, 
then $T_1\otimes\cdots\otimes T_m\ne0$ if and only if $T_k\ne 0$ for $k=1,\dots,m$, 
it is easily seen that condition~\eqref{fix_eq1}  
is equivalent to 
$$T_k:=\sum_{\sigma\in S_{q_k}} V_{0k}(\sigma)\otimes \overline{U_{0k}}(\sigma)\ne0
\text{\ for\ }k=1,\dots,m.
$$
To see that the above condition is satisfied, note that  
the irreducibility of the representation~$U_{0k}$ implies  
$$\frac{1}{q_k!}\sum_{\sigma\in S_{q_k}} U_{0k}(\sigma)=\1\text{ for }k=1,\dots,m. $$
Hence  for any nonzero vector $v\in \cH_{U_{0k}}$   
and 
for every nonzero symmetric tensor $w\in S^{q_k}(\cH_{\pi_{j_k^0}})\subseteq \cH_{\pi_{j_k^0}}^{\otimes q_k}$   
one has 
$$\begin{aligned}
0
&\ne w\otimes (q_k! v)=w \otimes \sum_{\sigma\in S_{q_k}} \overline{U_{0k}}(\sigma)v
=\sum_{\sigma\in S_{q_k}} V_{0k}(\sigma)w\otimes \overline{U_{0k}}(\sigma)v \\
&=\Bigl(\sum_{\sigma\in S_{q_k}} V_{0k}(\sigma)\otimes \overline{U_{0k}}(\sigma)\Bigr)(w\otimes v).
\end{aligned}$$
This completes the proof of  Assertion~\eqref{exact_item2}. 

(iii) follows by Theorem~\ref{fix}, and we are done.
\end{prf}

\begin{ex}\mlabel{compact}
Consider the $C^*$-algebra $\cA=\fS_\infty(\cH)$ of compact operators on some 
infinite dimensional complex  Hilbert space. 
Then it is known that $\widehat{\cA}=\{[\pi]\}$, 
where $\pi$ is the tautological representation of~$\cA$. 
Then, for any $n\ge 2$, in Theorem~\ref{exact} one has $m=1$, 
and it follows that $S^n(\cA)\,\hat{}\ $ is parametrized by $\widehat{S_n}$, 
just as the holomorphic Schur--Weyl representations.  
See Subsection~\ref{SWsubsect} below for more details. 
\end{ex}

\subsection{Links to Schur--Weyl theory for irreducible representations}\mlabel{SWsubsect}

Let $\cA$ be any $C^*$-algebra and recall from the introduction that its {\it unitary group} is 
\[ \U(\cA):=(\1+\cA)\cap \U(M(\cA)),\] 
where $M(\cA)$ is the multiplier $C^*$-algebra of $\cA$,  
and $M(\cA)=\cA$ if $\cA$ is unital. 
It follows from \cite{BN12} that to every irreducible $*$-representation 
$\pi\colon\cA\to B(\ell^2(J))$ 
there corresponds an infinite family of unitary 
irreducible representations $\pi_\lambda^\cA$ 
of $\U(\cA)$ obtained by Schur--Weyl theory (see \cite[Def. 4.1]{BN12}),   
where $\lambda$  belongs to the additive group $\cP_{J}\cong \Z^{(J)}$ 
of all finitely supported $\Z$-valued functions on $J$.  
Here a set $J=J_{[\pi]}$ is fixed for very equivalence class 
$[\pi]\in\widehat{\cA}$. 
We then define the \emph{Schur--Weyl map} of $\cA$ as 
\begin{equation}\label{SW_eq0}
\Psi_{\cA}\colon \bigsqcup_{[\pi]\in\widehat{\cA}} \cP_{J_{[\pi]}}/S_{(J_{[\pi]})}
\to \U(\cA)\,\hat{} ,\quad 
[\lambda]\mapsto[\pi_\lambda^\cA] 
\quad\text{ for }[\lambda]\in\cP_{J_{[\pi]}}/S_{(J_{[\pi]})}\text{ and }[\pi]\in\widehat{\cA}, 
\end{equation}
where the notation is as in \cite[Th.~3.2]{BN12}.  
In particular, for $\lambda\in\cP_J$, we denote by $[\lambda]$ its equivalence class modulo 
the natural action of the group $S_{(J)}$ of finitely supported permutations of $J$.  
As an application of Theorem~\ref{factbij}  
(via Theorem~\ref{fix}) we will show in Theorem~\ref{SW-main} that 
\emph{the Schur--Weyl map $\Psi_{\cA}$ is always injective} but, if $\cA$ is separable 
and non-isomorphic to $\cK(\ell^2(\N))$, then $\hat\cA$ contains 
at least two elements and this implies that $\Psi_{\cA}$ is not surjective. 

On the other hand, using \cite[Th. 6.8]{BN15}, one easily sees that  
if $\lambda_j \geq 0$ for every $j$ and $n = \sum_j \lambda_j$, then 
the representation $\pi_\lambda^\cA$ uniquely extends to an 
irreducible $*$-representation of $S^n(\cA)$, 
denoted again by $\pi_\lambda^\cA$ and termed a \emph{Schur--Weyl representation} in the following, 
where by extension we actually mean a factorization 
through the $*$-morphism 
\[ \U(\cA)\hookrightarrow S^n(\cA),\quad  a\mapsto a^{\otimes n}.\] 
To begin with, we will show in Proposition~\ref{SW} 
how the equivalence classes $[\pi_\lambda^\cA]\in S^n(\cA)\,\hat{}\ $  
can be recovered in the setting of the above Theorem~\ref{factbij}.  
To this end we need the following general lemma. 

\begin{lemma}\mlabel{quasi-SW}
Let $(\pi,U)$ be a covariant representation of a $C^*$-dynamical 
system $(\cB,G,\alpha)$, with the finite group $G$, on a Hilbert space 
$\cH$  for which  $\pi$ is irreducible.  
For any unitary irreducible representation $\lambda\colon G\to B(\cV)$ define  
the Hilbert space $\cH^\cV:=B_2(\cV,\cH)\simeq\cH\otimes\cV^*$ and the $*$-representation 
$\pi^{\cV}\colon\cB\to B(\cH^\cV)$, $\pi^{\cV}(x)T:=\pi(x)T$. 
Then the following assertions hold: 
\begin{enumerate}[\rm(i)] 
\item\label{quasi-SW_item1} The pair $(\pi^{\cV},U\otimes\overline{\lambda})$ 
is an irreducible covariant representation of $(\cB,G,\alpha)$.   
\item\label{quasi-SW_item2} Let $\iota\colon\cB^G\to\cB \subeq \cB\rtimes G$ 
be the canonical embedding as in {\rm Proposition~\ref{cross}\eqref{cross_item1}}.
If $\Hom_G(\cV,\cH)\ne\{0\}$, 
then this is the essential space of  
$(\pi^{\cV}\rtimes(U\otimes\overline{\lambda}))\circ\iota\colon\cB^G\to B(\cH^\cV)$ 
and the cor\-res\-ponding nondegenerate representation 
of $\cB^G$ on $\Hom_G(\cV,\cH)$ is irreducible 
and a subrepresentation of $\pi^{\cV}\vert_{\cB^G}$.  
\end{enumerate}
\end{lemma}

\begin{prf}
For every $g\in G$ and $x\in\cB$ one has 
$$(U\otimes\overline{\lambda})(g)\pi^{\cV}(x)(U\otimes\overline{\lambda})(g^{-1}) =(U_g\pi(x)U_{g^{-1}})\otimes\id_{\cV^*}=\pi(\alpha_g(x))\otimes\id_{\cV^*}
=\pi^{\cV}(\alpha_g(x)),  $$
so that $(\pi^{\cV},U\otimes\overline{\lambda})$ is a covariant representation of $(\cB,G,\alpha)$. 

To prove that $(\pi^{\cV},U\otimes\overline{\lambda})$ is irreducible, we must check that 
$\pi^{\cV}(\cB)'\cap(U\otimes\overline{\lambda})_G'=\C\id_{\cH\otimes\cV^*}$. 
One has $\pi^{\cV}(\cB)'=(\pi(\cB)\otimes\C \id_{\cV^*})'=\pi(\cB)'\otimes B(\cV^*) =\C\id_{\cH}\otimes B(\cV^*)$ 
because $\pi$ is irreducible. 
Hence we must prove that if $T\in B(\cV^*)$ satisfies $\id_{\cH}\otimes T\in (U\otimes\overline{\lambda})_G'$, 
then $T\in\C\id_{\cV^*}$. 
The condition on $T$ is equivalent to $T\in\overline{\lambda}_G'$ hence, using the hypothesis that $\lambda$ is an irreducible representation 
hence so is $\overline{\lambda}$, 
one obtains $T\in\C\id_{\cV^*}$, and this completes the proof of the first assertion. 

For the second assertion, it follows by Proposition~\ref{cross}(ii) and (iii) 
that $\Hom_G(\cV,\cH)$, 
that is, $(\cH\otimes\cV^*)^G$,  
is the essential space of $(\pi^{\cV}\rtimes(U\otimes\overline{\lambda}))\circ\iota$, 
and the corresponding nondegenerate representation of $\cB^G$ is irreducible because 
$\pi^{\cV}\rtimes(U\otimes\overline{\lambda})$ is irreducible. 
To see that this irreducible representation of $\cB^G$ is a subrepresentation of $\pi^\cV\vert_{\cB^G}$, 
that is, it is equal to $x\mapsto\pi^\cV(x)\vert_{\Hom_G(\cV,\cH)}$, 
note that, for $x\in\cB^G$,  one has 
$$\begin{aligned}
((\pi^{\cV}\rtimes(U\otimes\overline{\lambda}))\circ\iota)(x)
&=\frac{1}{\vert G\vert}\sum_{g\in G}(\pi^{\cV}\bigl(\iota(x)(g)\bigr) (U\otimes\overline{\lambda})_g
=\frac{1}{\vert G\vert}\sum_{g\in G}\pi^{\cV}(x) (U\otimes\overline{\lambda})_g \\
&=\pi^{\cV}(x)P=P\pi^{\cV}(x), 
\end{aligned}$$
where $P:=\frac{1}{\vert G\vert}\sum\limits_{g\in G}(U\otimes\overline{\lambda})_g$ 
is the orthogonal projection of $\cH\otimes\cV^*$ onto $(\cH\otimes\cV^*)^G$. 
This completes the proof. 
\end{prf}

\begin{rem}
In connection with Lemma~\ref{quasi-SW}, we note that the factor covariant representations $(\pi,U)$ 
of $(\cB,G,\alpha)$ for which the irreducible decomposition of $U$ contains a fixed irreducible representation $\lambda$ of $G$ 
were classified in \cite[Th. 1]{La80} even for compact groups~$G$. 
\end{rem}

\begin{prop}\mlabel{SW}
Let $\cA$ be a unital $C^*$-algebra and fix some integer $n\ge1$ and  
a unitary irreducible representation $\lambda\colon S_n\to B(\cV)$. 
For any irreducible $*$-representation $\pi\colon\cA\to B(\cH)$ 
denote by $W(\pi)\colon S_n\to B(\cH^{\otimes n})$ the corresponding permutation representation 
and define  
\[ (\pi^{\otimes n})^\cV\colon\cA^{\otimes n}\to B(\cH^{\otimes n}\otimes\cV^*), 
\quad (\pi^{\otimes n})^\cV(x):=\pi^{\otimes n}(x)\otimes\id_{\cV^*}.\]
Then $((\pi^{\otimes n})^{\cV},W(\pi)\otimes\overline{\lambda})$ 
is an irreducible covariant representation of $(\cA^{\otimes n},S_n,\alpha)$ 
and the image of its $S_n$-orbit through the map $\cI$ of {\rm Theorem~\ref{fix}} is 
$[\pi_\lambda^\cA]\in S^n(\cA)\,\hat{}\ \hookrightarrow\U(\cA)\,\hat{} $. 
\end{prop}

\begin{prf}
The covariant representation $((\pi^{\otimes n})^{\cV},W(\pi)\otimes\overline{\lambda})$ is irreducible by Lemma~\ref{quasi-SW}, 
because $\pi^{\otimes n}\colon\cA^{\otimes n}\to B(\cH^{\otimes n})$ is irreducible. 
Using the notation of Definition~\ref{masamichi}, 
one has $\pi_0=\pi^{\otimes n}$ and $G=G_{\pi_0}=S_n$, hence $\Ind_{G_{\pi_0}}^G(\pi_0)=\pi_0$. 
Moreover, $\omega_0\equiv 1$, $V_0=W(\pi)$, 
and $U_0=\lambda$. 
Condition~\eqref{fix_eq1} is satisfied, since this can be proved just as the relation $T_k\ne0$ in 
the last part of the proof of Theorem~\ref{exact}\eqref{exact_item2}. 
The image of the $S_n$-orbit of the equivalence class of $(\pi^{\cV},W(\pi)\otimes\overline{\lambda})$ under $\cI$
is the equivalence class of the nondegenerate part of the representation
$(\pi_0\rtimes (V_0\otimes\overline{U_0}))\circ\iota=(\pi^{\otimes n}\rtimes (W(\pi)\otimes\overline{\lambda}))\circ\iota$. 
This nondegenerate part is unitarily equivalent to the representation 
\begin{equation}\label{SW_proof_eq1}
S^n(\cA)=(\cA^{\otimes n})^{S_n}\to B(\Hom_{S_n}(\cV,\cH^{\otimes n})), \quad
x\mapsto \pi^{\otimes n}(x)\otimes\id_{\cV^*}
\end{equation}
by Lemma~\ref{quasi-SW}\eqref{quasi-SW_item2}, 
since 
$$(\pi^{\otimes n}(x)\otimes\id_{\cV^*})(T)=\pi^{\otimes n}(x)T\text{ for every }
T\in\Hom_{S_n}(\cV,\cH^{\otimes n})\simeq (\cH^{\otimes n}\otimes\cV^*)^{S_n}.$$ 
But then, by composing \eqref{SW_proof_eq1} with the power map 
$\U(\cA)\to S^n(\cA), a\mapsto a^{\otimes n}$ 
one obtains a unitary irreducible representation of $\U(\cA)$ 
(compare \cite[Cor. 6.10]{BN15})
on the Hilbert space 
$\Hom_{S_n}(\cV,\cH^{\otimes n})=:\bS_\lambda(\cH^{\otimes n})$  
which is unitarily equivalent to the Schur--Weyl representation $\pi_\lambda^\cA$  
(see \cite[Rem. A.6]{BN12}). 
\end{prf}

In Assertion~(\ref{SW-main_item2}) of the following statement we use the notation 
$\cC_1(\cA^{\otimes n},S_n,\alpha)$ as introduced in Example~\ref{irregular}. 

\begin{thm}\mlabel{SW-main}
For any $C^*$-algebra $\cA$, the following assertions hold: 
\begin{enumerate}[\rm(i)]
\item\label{SW-main_item1} The Schur--Weyl map $\Psi_{\cA}$  from 
\eqref{SW_eq0} is injective. 
\item\label{SW-main_item2} One has 
$\Im\Psi_{\cA}\subseteq\bigcup\limits_{n\ge 0}S^n(\cA)\,\hat{}\ \hookrightarrow\U(\cA)\,\hat{} $
and $\Im\Psi_{\cA}\cap S^n(\cA)\,\hat{}\ $ consists of 
restrictions of representations which belong to $\cC_1(\cA^{\otimes n},S_n,\alpha)$, for all $n\ge1$. 
\item The map $\Psi_{\cA}$ is bijective if $\cA$ is isomorphic to 
the $C^*$-algebra of compact operators on some complex Hilbert space.  
\item If $\pi_j\colon \cA\to B(\cH_j), j=1,2,$ are 
inequivalent irreducible $*$-representations, then the map 
$a\mapsto \pi_1(a)\otimes\pi_2(a)$ 
is a unitary irreducible representation of $\U(\cA)$ on $\cH_1\otimes\cH_2$ 
whose equivalence class does not belong to the image of the Schur--Weyl map. 
\end{enumerate}
\end{thm}

\begin{prf} (i), (ii)  
Recall from \eqref{SW_eq0} that the Schur--Weyl map is 
$$\Psi_{\cA}\colon \bigsqcup_{{[\pi]}\in\widehat{\cA}} \cP_{J_{[\pi]}}/S_{(J_{[\pi]})}
\to \U(\cA)\,\hat{} ,\quad 
[\lambda]\mapsto[\pi_\lambda^\cA] 
\quad\text{ for }[\lambda]\in\cP_{J_{[\pi]}}/S_{(J_{[\pi]})}\text{ and }[\pi]\in\widehat{\cA}.$$
To prove that $\Psi_\cA$ is injective, 
let $[\pi_j]\in \widehat{\cA}$ and $\lambda_j\in\cP_{J_{[\pi_j]}}$ for $j=1,2$ with 
$[(\pi_1)_{\lambda_1}^\cA]=[(\pi_2)_{\lambda_2}^\cA]\in \U(\cA)\,\hat{}$. 
We must check that $[\pi_1]=[\pi_2] :=[\pi]\in \widehat{\cA}$ 
and $[\lambda_1]=[\lambda_2]\in \cP_{J_{[\pi]}}/S_{(J_{[\pi]})}$. 
 
Denote also by $\lambda_j\colon S_n\to B(\cV_j)$ the unitary irreducible representation associated with 
$\lambda_j\in\cP_{J_{[\pi]}}$, where $n\ge \max\{\vert\supp\,\lambda_1\vert,\vert\supp\,\lambda_2\vert\}$.  
It follows by Proposition~\ref{SW} that $[\pi_\lambda]\in S^n(\cA)\,\hat{}\ $ is 
the image through the bijection $\cI$ (Theorem~\ref{fix}) 
of the irreducible covariant representation 
$((\pi_j^{\otimes n})^{\cV},W(\pi_j)\otimes\overline{\lambda_j})$ 
 of $(\cA^{\otimes n},S_n,\alpha)$. 
Then $[(\pi_1)_{\lambda_1}^\cA]=[(\pi_2)_{\lambda_2}^\cA]\in S^n(\cA)\,\hat{}\ \hookrightarrow \widehat{\cA}$, 
implies 
\[ [(\pi_1^{\otimes n})^{\cV_1}\rtimes (W(\pi_1)\otimes\overline{\lambda_1})]
=[(\pi_2^{\otimes n})^{\cV_2}\rtimes (W(\pi_2)\otimes\overline{\lambda_2})]
\in (\cA^{\otimes n}\rtimes S_n)\,\hat{}\ .\]
Now we can use Example~\ref{irregular} and the injectivity of the map \eqref{new_prf0} in Lemma~\ref{lprf0} 
to obtain $[\pi_1]=[\pi_2]$, that is, $[\pi_1] = [\pi_2] :=[\pi]\in \widehat{\cA}$. 
It then follows by \cite[Th. 4.4]{BN12} that also $[\lambda_1]=[\lambda_2]\in \cP_{J_{[\pi]}}/S_{(J_{[\pi]})}$, 
and this completes the proof of the fact that $\Psi_\cA$ is injective. 

(iii) If $\cA$ is isomorphic to 
the $C^*$-algebra of compact operators on some complex Hilbert space, then $\Psi_\cA$ is surjective by 
\cite[Th. 3.21]{Ne14}.  

(iv) Now let $\pi_j\colon \cA\to B(\cH_j)$ for $j=1,2$ be inequivalent irreducible $*$-representations 
and define the irreducible $*$-representation $\pi_0:=\pi_1\otimes\pi_2\colon \cA^{\otimes 2}\to B(\cH_1\otimes\cH_2)$. 
Using again the injectivity of the map \eqref{new_prf0} in Lemma~\ref{lprf0} it follows
 that the isotropy group of $[\pi_0]$ for the natural action of $S_2$ on $\widehat{\cA^{\otimes 2}}$ is $\{\1\}$, 
 that is, $[\pi]$ belongs to the class $\cC_0(\cA^{\otimes 2},S_2,\alpha)$ from Example~\ref{regular}. 
This shows that the image of $[\pi_0]$ under $\cI$ does not belong to the image of the Schur--Weyl map $\Psi_\cA$. 
It only remains to compute the image of $[\pi_0]$ under the map $\cI$. 
To this end we first compute  
the corresponding induced representation $(\pi,U)$ of $(\cA^{\otimes 2},S_2,\alpha)$ 
given by Definition~\ref{Def3.1}. 
Since $G_0=\{\1\}$, we have $\cH=\cH_0\oplus\cH_0$ thought of as column vectors and 
$$U\colon S_2\simeq\Z_2\to B(\cH) \text{ with }U_0=\1\text{ and }
U_1=\begin{pmatrix} 0 & 1 \\ 
                     1 & 0
     \end{pmatrix}$$
and also 
$$\pi\colon\cA^{\otimes 2}\to B(\cH),\quad 
\pi(a_1\otimes a_2)=\begin{pmatrix} \pi_1(a_1)\otimes\pi_2(a_2) & 0 \\ 
                     0 & \pi_1(a_2)\otimes\pi_2(a_1)
     \end{pmatrix}$$ 
for all $a_1,a_2\in\cA$. 
If we view the elements of the crossed product $\cA^{\otimes 2}\rtimes S_2$ as pairs of elements of $\cA^{\otimes 2}$ 
with the multiplication given by 
\[ \begin{aligned}
(a_1\otimes a_2,  & b_1\otimes b_2)\cdot (c_1\otimes c_2, d_1\otimes d_2) \\
& =((a_1c_1)\otimes(a_2c_2)+(b_1d_1)\otimes (b_2d_2),(a_1d_1)\otimes(a_2d_2)+(b_1c_2)\otimes (b_2c_1)) 
\end{aligned}\]
then the irreducible representation $\pi\rtimes U\colon \cA^{\otimes 2}\rtimes S_2\to B(\cH)$ 
is given by 
$$\begin{aligned}
(\pi\rtimes U)(a_1\otimes a_2, b_1\otimes b_2)
& =\pi(a_1\otimes a_2)U_0+\pi(b_1\otimes b_2)U_1 \\
&=
\begin{pmatrix} \pi_1(a_1)\otimes\pi_2(a_2) & \pi_1(b_2)\otimes\pi_2(b_1) \\ 
                \pi_1(b_1)\otimes\pi_2(b_2) & \pi_1(a_2)\otimes\pi_2(a_1)
\end{pmatrix}
\end{aligned}$$
and by composing $\pi\rtimes U$ with 
$\iota\colon S^2(\cA)\to \cA^{\otimes 2}\rtimes S_2$, $\iota(a\otimes a)=(a\otimes a, a\otimes a)$, 
we obtain a degenerate representation on $\cH$ whose essential subspace 
is the image of $\begin{pmatrix}
\1 & \1 \\
\1 & \1
\end{pmatrix}$, i.e., the ``diagonal subspace'' of $\cH=\cH_0\oplus\cH_0$. 
The nondegenerate part of $(\pi\rtimes U)\circ\iota$ is unitarily equivalent to the representation 
$$S^2(\cA)\to B(\cH_0),\quad a\otimes a\mapsto \pi_1(a)\otimes\pi_2(a)$$
(and is irreducible by Proposition~\ref{cross}\eqref{cross_item3}). 
Composing it further with the embedding $\U(\cA)\to S^2(\cA)$, $a\mapsto a\otimes a$, 
we obtain the irreducible representation of $\U(\cA)$ from the statement. 
This concludes the proof. 
 \end{prf}

\section{Schur--Weyl theory for factor representations} 
\mlabel{sec:6}

In Propositions \ref{fact1}--\ref{fact2} below we give a generalization of 
\cite[Th.~1(1)]{Nes13} 
dealing with factors of type ${\rm II}_1$ to more general von Neumann algebras, 
and this also provides a partial generalization of \cite[Th. 4.1]{EnIz15} 
where a Schur--Weyl duality property was established for the standard representation of the hyperfinite factor of type ${\rm II}_1$.  
In the case $\cM=B(\cH)$ the following proposition also recovers the Schur--Weyl duality studied in \cite{BN12}. 

\begin{prop}\mlabel{fact1}
Let $\cM\subseteq B(\cH)$ be a von Neumann algebra and for any $n\ge 1$ consider 
the unitary permutation representation 
$V\colon S_n\to B(\cH^{\otimes n})$ and its corresponding 
action by $*$-automorphisms $\alpha\colon S_n\to\Aut(\cM^{\otimes n})$. 
Then for the homomorphism of multiplicative $*$-semigroups  
$\Gamma\colon\cM\to (\cM^{\otimes n})^{S_n}$, $\Gamma(a):=a^{\otimes n}$,  
one has 
$$\Gamma(\U(\cM))''=(\cM^{\otimes n})^{S_n}\quad\text{and}\quad
\Gamma(\U(\cM))'=(V_{S_n}\cup(\cM')^{\otimes n})''.$$ 
\end{prop}

\begin{prf}
The inclusion $\Gamma(\U(\cM))''\subseteq(\cM^{\otimes n})^{S_n}$ is clear. 
For the opposite inclusion we use the differential of $\Gamma$ at $0\in\cM$, 
$$\dd\Gamma\colon \cM\to (\cM^{\otimes n})^{S_n},\quad 
\dd\Gamma(a)=\sum_{k=1}^n\1^{\otimes (k-1)}\otimes a\otimes\1^{\otimes (n-k)}.$$
By the formulas $\dd\Gamma(a)=\frac{1}{i}\derat0\Gamma(e^{it a})$ and $\Gamma(e^{ia})=e^{i\dd\Gamma(a)}$ for all $a=a^*\in\cM$, 
one easily obtains the equality 
$\Gamma(\U(\cM))''=\dd\Gamma(\cM)''$, 
hence it suffices to prove that $(\cM^{\otimes n})^{S_n}\subseteq \dd\Gamma(\cM)''$. 
As one has the surjective map 
$$E\colon\cM^{\otimes n}\to (\cM^{\otimes n})^{S_n},\quad 
E(b)=\frac{1}{n!}\sum_{\sigma\in S_n}\alpha_\sigma(b)$$
we will have to prove
\begin{equation}\label{fact1_proof_eq1}
(\forall b_1,\dots,b_n\in\cM)\quad E(b_1\otimes\cdots\otimes b_n)\in \dd\Gamma(\cM)''. 
\end{equation}
We will prove this by recurrence after $r:=\vert\{j\in\{1,\dots,n\}\mid b_j\ne\1\}\vert$. 
If $r=1$, this is clear since $E(b_1\otimes\cdots\otimes b_n)=\frac{1}{n!}\dd\Gamma(b_{j_1})$, where $j_1$ is the unique 
$j\in\{1,\dots,n\}$ with $b_j\ne\1$. 
Now assume $r+1\le n$ and define $\ell_0:=\max\{j\in\{1,\dots,n\}\mid b_j\ne\1\}$, 
hence $b_{\ell_0}\ne\1$. 
Also define 
$c_j:=b_j$ if $j\in\{1,\dots,n\}\setminus\{\ell_0\}$ 
and $c_{\ell_0}:=\1$. 
Then one has 
$$\begin{aligned}
\dd\Gamma(b_{\ell_0})\cdot n!E(c_1\otimes\cdots\otimes c_n)
&=\Bigl(\sum_{k=1}^n \1^{\otimes (k-1)}\otimes b_{\ell_0}\otimes\1^{\otimes (n-k)}\Bigr) 
\Bigl(\sum_{\sigma\in S_n} c_{\sigma(1)}\otimes\cdots\otimes c_{\sigma(n)}\Bigr) \\
&=\sum_{k=1}^n
\sum_{\sigma\in S_n} c_{\sigma(1)}\otimes\cdots\otimes(b_{\ell_0} c_{\sigma(k)})\otimes\cdots \otimes c_{\sigma(n)}
\end{aligned}
$$
If we split the above sum according to the pairs $(\sigma,k)\in S_n\times\{1,\dots,n\}$ 
which satisfy 
$\sigma(k)=\ell_0$ and $\sigma(k)\ne\ell_0$,  
respectively, 
then we obtain 
$$\begin{aligned}
\dd\Gamma(b_{\ell_0})\cdot n!E(c_1\otimes\cdots\otimes c_n)
=& n!E(b_1\otimes\cdots \otimes b_n) \\
&+n!\sum_{k\in\{1,\dots,n-1\}\setminus\{\ell_0\}}E(b_1\otimes\cdots\otimes(b_{\ell_0}b_k)\otimes\cdots\otimes b_n) 
\end{aligned}
$$
where in every summand of the second sum, the $\ell_0$-th factor in the tensor product is equal to~$\1$. 
Solving the above equation for $E(b_1\otimes\cdots \otimes b_n)$ and using the recurrence hypothesis, 
one completes the proof of \eqref{fact1_proof_eq1}, and we have already seen above that this implies the assertion. 

It remains to prove the second equality from the statement. 
Using the commutation theorem for von Neumann algebras $(\cM_1\otimes\cM_2)'=\cM_1'\otimes\cM_2'$ 
and the Bicommutant Theorem, one obtains  
$$(V_{S_n}\cup(\cM')^{\otimes n})'
=(V_{S_n}\cup(\cM^{\otimes n})')'
=V_{S_n}'\cap\cM^{\otimes n}
=(\cM^{\otimes n})^{S_n}
=\Gamma(\U(\cM))''$$
where we used also the first of the asserted equalities from the statement, which was already proved. 
Now the second equality from the statement follows by taking commutants in the above equalities and using once again the Bicommutant Theorem. 
\end{prf}

The following proposition shows that, beyond the classical Schur--Weyl theory when $\cM=B(\cH)$, 
the picture is completely different.  

\begin{prop}\mlabel{fact2}
Assume the setting of \emph{Proposition~\ref{fact1}}. 
If $\cM$ is a factor of type II$_1$,  II$_\infty$, or III, then $\Gamma\colon \U(\cM)\to B(\cH^{\otimes n})$, $\Gamma(u)=u^{\otimes n}$,  
is a factor representation of the same type as $\cM$.  
\end{prop}

\begin{prf} 
By Proposition~\ref{fact1} it suffices to prove that 
 $(\cM^{\otimes n})^{S_n}$ is a factor of the same type as~$\cM$. 
 
Since $\cM$ is not a factor of type~I, it follows by \cite[Th. 5]{Sa75} that for every $\sigma\in S_n\setminus\{\1\}$ 
the automorphism $\alpha(\sigma)$ of $\cM^{\otimes n}$ is not inner, 
that is, it is not of the form $x\mapsto uxu^*$ for any unitary element $u\in\cM$. 
Then, since $\cM$ is a factor, 
 it follows by \cite[Prop. 2.6.7]{Sa71} that $\cM^{\otimes n}$ is a factor, 
hence it follows by \cite[Cor. to Prop. II.3]{Au76} that $(\cM^{\otimes n})^{S_n}$ is in turn a factor. 

To check that $(\cM^{\otimes n})^{S_n}$ has the same type as~$\cM$, 
we discuss below separately the cases that can occur, 
using the faithful normal conditional expectation 
$$E\colon \cM^{\otimes n}\to (\cM^{\otimes n})^{S_n},\quad E(a)=\frac{1}{n!}\sum_{\sigma\in S_n}\alpha_\sigma(a).$$
If $0\le a\in \cM^{\otimes n}$ then, using the fact that $0\le \alpha_\sigma(a)$ if $\sigma\in S_n\setminus\{\1\}$ 
and $a=\alpha_\sigma(a)$ if $\sigma=\1$, one obtains $a\le n! E(a)$, 
hence the conditional expectation $E$ has finite index in the sense of \cite{FrKi98}. 
It then follows by \cite[Prop. 2.2]{FrKi98} that the factors $\cM^{\otimes n}$ and $(\cM^{\otimes n})^{S_n}$ 
have the same type II$_1$, II$_\infty$, or III. 

If $\cM$ is of type III, then $\cM^{\otimes n}$ is of type III by \cite[Th. 2.6.4]{Sa71}. 
If $\cM$ is of type II$_1$ or~II$_\infty$, then $\cM^{\otimes n}$ is of the same type by \cite[Prop. 2.6.1/3]{Sa71}. 
This completes the proof. 
\end{prf}

\begin{rem}
With the notation of the proof of Proposition~\ref{fact2},  
the index of the conditional expectation $E$ is $\Ind E=\vert S_n\vert=n!$ by \cite[Ex. 1.2]{Lo92}.  
Moreover, if $\cM$ is of type III$_\lambda$ for some $\lambda\in[0,1]$, then 
more refined information on the type of $(\cM^{\otimes n})^{S_n}$ can be obtained 
by \cite[Th. 2.7]{Lo92}. 
\end{rem}

The following theorem extends \cite[Th. 1(2)]{Nes13} to infinite factors. 

\begin{thm}\mlabel{fact3}
In the setting of \emph{Proposition~\ref{fact2}} we define 
for every non-zero projection $p=p^2=p^*\in\Gamma(\U(\cM))'$ with range $\cH_p= p \cH$ 
the representation 
\[ {\Gamma_p\colon \U(\cM)\to \U(\cH_p), \quad \Gamma_p(u):=\Gamma(u)\vert_{\cH_p}}.\]
\begin{enumerate}[\rm(i)]
\item $\Gamma_p$ is a factor representation of the same type II or III as $\cM$ and~$\Gamma_p\approx\Gamma$. 
\item If $\cM$ is a factor of type II, then for any faithful normal trace $\tau$ on $\Gamma(\U(\cM))'$   
and finite projections $p_j=p_j^2=p_j^*\in\Gamma(\U(\cM))'$ for $j=1,2$, 
one has $\Gamma_{p_1}\le \Gamma_{p_2}$ if and only if $\tau(p_1)\le\tau(p_2)$, 
and $\Gamma_{p_1}\simeq \Gamma_{p_2}$ if and only if $\tau(p_1)=\tau(p_2)$. 
\item If $\cM$ is a countably decomposable factor of type III, 
then $\Gamma_p\simeq\Gamma$. 
\end{enumerate}
\end{thm}

\begin{prf} (i) 
It is easily seen that $\Gamma_p(\U(\cM))''=\Gamma(\U(\cM))_p$ 
(reduced von Neumann algebra), 
hence  $\Gamma_p$ is a factor representation, by \cite[Th. 3.13]{SZ79}. 
Moreover, the reduction map $\Gamma(\U(\cM))\to \Gamma(\U(\cM))_p$ is a $*$-isomorphism by \cite[Prop. 3.14]{SZ79}, 
hence $\Gamma$ and $\Gamma_p$ are quasi-equivalent factor representations. 

(ii) If $\cM$ is a factor of type II, then $\Gamma(\U(\cM))''$ is a factor of type II by Proposition~\ref{fact2}, 
hence by \cite[Ch. I, \S 6, no. 8, Cor. 1 of Prop. 13]{Dix69} also $\Gamma(\U(\cM))'$ is a factor of type II, 
and then it has a faithful normal trace $\tau$. 
If $p_j=p_j^2=p_j^*\in\Gamma(\U(\cM))'$ are finite projections for $j=1,2$, then one has $\tau(p_1)\le\tau(p_2)$ 
(respectively $\tau(p_1)=\tau(p_2)$) if and only if 
$p_1\prec p_2$ (respectively $p_1\sim p_2$) in $\Gamma(\U(\cM))'$ by \cite[Ch. III, \S 2, Prop. 13]{Dix69}, 
which is further equivalent to $\Gamma_{p_1}\le \Gamma_{p_2}$ by \cite[Cor. 5.1.4]{Dix64} 
(respectively $\Gamma_{p_1}=\Gamma_{p_2}$ by \cite[Cor. 5.1.3]{Dix64}). 

(iii) If $\cM$ is a countably decomposable factor of type III 
and $0\ne p_j=p_j^2=p_j^*\in\Gamma(\U(\cM))'$ for $j=1,2$, 
then $p_1\sim p_2$ in $\Gamma(\U(\cM))'$ by \cite[Prop. 2.2.14]{Sa71}, hence $\Gamma_{p_1}$ and $\Gamma_{p_2}$ are unitarily  equivalent again by \cite[Cor. 5.1.3]{Dix64}. 
\end{prf}

\begin{rem}\mlabel{fact4}
It follows by Proposition~\ref{fact1} that to every central projection $\chi\in\C[S_n]$ there corresponds an orthogonal projection 
$V(\chi)\in\Gamma(\U(\cM))'$. 
In the case when $\cM$ is a factor of type II, the traces of these projections were computed in the proof of \cite[Th. 1(3)]{Nes13} and in \cite[Cor. 5]{DaKa06}. 

When $\cM$ is a factor of type II$_1$ with its faithful normal tracial state~$\tau$, 
the character of the representation $\Gamma$ is $(\tau\vert_{\U(\cM)})^n$ (which belongs to the list in \cite[Th. 1.5]{EnIz15}) 
and this is equal to the character of the representation $\Gamma_p$ if $0\ne p=p^2=p^*
\in\Gamma(\U(\cM))'$, 
because, by Theorem~\ref{fact3}, the representations $\Gamma$ and $\Gamma_p$ are quasi-equivalent. 
\end{rem}

\subsection*{Acknowledgment} 
We wish to thank Jan Stochel for informing us about his paper \cite{St92}. 

\end{document}